\documentclass[
               DIV=15,
               ]{scrartcl}  

\usepackage{xcolor} 
\usepackage{psfrag} 
\usepackage{subcaption}
\usepackage{url}
\usepackage{amsmath}
\usepackage{graphics}
 \usepackage{pstool} 
\usepackage{ulem}

\usepackage[colorinlistoftodos,prependcaption,textsize=tiny,backgroundcolor=blue!10!white]{todonotes}

\usepackage[square,sort,comma,numbers]{natbib}
 
\usepackage[utf8]{inputenc}

\usepackage{pgfplots}

\usepackage{tikz}
\usepackage{tikzscale}

\usepackage{algorithm}
\usepackage{algpseudocode}

\usetikzlibrary{spy}
\usetikzlibrary{arrows,shapes}    
\usetikzlibrary{backgrounds}  
\usetikzlibrary{decorations}
\usetikzlibrary{positioning}
\usetikzlibrary{patterns}
\usetikzlibrary{calc} 

\usepackage{booktabs}   
\usepackage{makecell} 
\usepackage{colortbl}   

\usepackage{xspace}
\usepackage{color}     
\usepackage{setspace}   
 
\usepackage{soul}
\usepackage{ulem}
\usepackage{currfile}
\usepackage{import}

\usepackage{authoraftertitle} 
\usepackage{authblk} 

\usepackage{cleveref}

\usepackage[titletoc,toc,title]{appendix} 

\usepackage{courier}

\pgfplotsset{compat=1.8}

\newcommand{\findmax}[3]{
    \pgfplotstablesort[sort key={#2},sort cmp={float >}]{\sorted}{#1}%
    \pgfplotstablegetelem{0}{#2}\of{\sorted}%
    \let #3=\pgfplotsretval%
}

\pgfplotscreateplotcyclelist{myCycleListBlackAndWhite}
{%
  {black, mark=o},
  {black, mark=square},
  {black, mark=triangle},
  {black, mark=diamond}, 
  {black, mark=pentagon}, 
  {black, mark=*},
  {black, mark=square*},
  {black, mark=triangle*},
  {black, mark=diamond*}, 
  {black, mark=pentagon*}, 
}

\definecolor{darkgreen}{rgb}{0,0.4,0} 
\definecolor{darkbrown}{rgb}{0.5, 0.396, 0.09}

\definecolor{c1}{rgb}{0.0, 0.4196078431372549, 0.6431372549019608}
\definecolor{c2}{rgb}{1.0, 0.5019607843137255, 0.054901960784313725}
\definecolor{c3}{rgb}{0.6705882352941176, 0.6705882352941176,
0.6705882352941176} \definecolor{c}{rgb}{0.34901960784313724, 0.34901960784313724, 0.34901960784313724}
\definecolor{c4}{rgb}{0.37254901960784315, 0.6196078431372549,
0.8196078431372549} \definecolor{c}{rgb}{0.7843137254901961, 0.3215686274509804, 0.0}
\definecolor{c5}{rgb}{0.5372549019607843, 0.5372549019607843,
0.5372549019607843} \definecolor{c}{rgb}{0.6352941176470588, 0.7843137254901961, 0.9254901960784314}
\definecolor{c6}{rgb}{1.0, 0.7372549019607844, 0.4745098039215686}
\definecolor{c7}{rgb}{0.8117647058823529, 0.8117647058823529,
0.8117647058823529}

\pgfplotscreateplotcyclelist{myCycleListColor}
{%
  {blue, mark=o},
  {red, mark=square},
  {darkbrown, mark=triangle},
  {darkgreen, mark=diamond},
  {black, mark=pentagon},
  {blue, mark=*},
  {red, mark=square*},
  {darkbrown, mark=triangle*},
  {darkgreen, mark=diamond*},
  {black, mark=pentagon*},
}

\pgfplotsset{every axis/.append style= 
              {
                font=\small,
                mark size=2,
                line width = 0.1,
                legend style={font=\small, mark size=3, draw=none, fill=none},
                legend cell align=left,
                cycle list name=myCycleListColor,
              }
            }
            
\pgfplotstableset
{
  every odd row/.style={before row={\rowcolor[gray]{0.9}}},
  every head row/.style=
  {
    before row=
    {%
      \toprule
    },
    after row=\midrule
  },
  every last row/.style=
  {
    after row=\bottomrule
  }
}

\newif\ifdrawboundingbox

\usetikzlibrary{external} 
\tikzset{external/system call={pdflatex \tikzexternalcheckshellescape
-halt-on-error -interaction=batchmode -jobname "\image" "\texsource"}} 

\pgfkeys
{ 
  /pgf/number format/.cd,
  fixed,
  set thousands separator={\,},
  1000 sep in fractionals,
}

\newcommand{\reva}[1]{\color{black}#1\color{black}}
\newcommand{\revb}[1]{\color{black}#1\color{black}}

\newcolumntype{C}[1]{>{\centering\arraybackslash}m{#1}}
\newcolumntype{R}[1]{>{\raggedright\arraybackslash}m{#1}}
\newcolumntype{L}[1]{>{\raggedleft\arraybackslash}m{#1}}

\newcommand{\delete}[1]{\xspace}

\title{Enforcing essential boundary conditions on domains defined by point clouds}

 \author[1]{Frank Hartmann \thanks{frank.hartmann@tum.de}}
 \author[2]{Stefan Kollmannsberger \thanks{stefan.kollmannsberger@tum.de, Corresponding author}}

\affil[1]{Technische Universit\"at M\"unchen,
  Arcisstr. 21, 80333 M\"unchen, Germany}

\affil[2]{Chair of Computational Modeling and Simulation,
 Technische Universit\"at M\"unchen,
  Arcisstr. 21, 80333 M\"unchen, Germany}

\drawboundingboxtrue
\drawboundingboxfalse 

\begin{document}  
\normalem
\maketitle  
  
\vspace{-2.5cm} 
\hrule 
\section*{Abstract}
This paper develops and investigates a new method for the application of Dirichlet boundary conditions for computational models defined by point clouds. Point cloud models often stem from laser or structured-light scanners which are used to scan existing mechanical structures for which CAD models either do not exist or from which the artifact under investigation deviates in shape or topology. Instead of reconstructing a CAD model from point clouds via surface reconstruction and a subsequent boundary conforming mesh generation, a direct analysis without pre-processing is possible using embedded domain finite element methods. These methods use non-boundary conforming meshes which calls for a weak enforcement of Dirichlet boundary conditions. For point cloud based models, Dirichlet boundary conditions are usually imposed using a diffuse interface approach. This leads to a significant computational overhead due to the necessary computation of domain integrals. Additionally, undesired side effects on the gradients of the solution arise which can only be controlled to some extent. This paper develops a new sharp interface approach for point cloud based models which avoids both issues. The computation of domain integrals is circumvented by an implicit approximation of corresponding Voronoi diagrams of higher order and the resulting sharp approximation avoids the side-effects of diffuse approaches. \revb{Benchmark examples from the graphics as well as the computational mechanics community are used to verify the algorithm. All algorithms are implemented in the FCMLab framework and provided at~\url{https://gitlab.lrz.de/cie_sam_public/fcmlab/}}. Further, we discuss challenges and limitations of point cloud based analysis w.r.t. application of Dirichlet boundary conditions.
 
\vspace{0.25cm}
\noindent \textit{Keywords:} finite cell method, point clouds, essential boundary conditions 
\vspace{0.25cm}
\hrule 
\newpage

\tableofcontents
\newpage


\section{Introduction} \label{sec:intro}

\subsection{Geometric Models}

Geometric models form the basis of any numerical analysis in computational mechanics. Popular geometric models in this field are based on Constructive Solid Geometry (CSG), Boundary representations (B-rep) or a mixture of both~\cite{Mantyla1988}. An in-depth discussion of the advantages and disadvantages w.r.t. computational mechanics is given, for example, in~\cite{Wassermann2019}. An interesting extension of B-rep models to volumes are V-rep models~\cite{Massarwi2016} which allow for the explicit description of the volume itself. To analyze the mechanical behavior of these structures, it is common engineering practice to generate them in Computer Aided Design systems prior to the production of an artifact.  

By contrast, for already existing structures it is often the case that a corresponding CAD model does not exist or that the produced artifact deviates from the CAD geometry. One way of obtaining a simple computer based volumetric model of an existing structure is to scan it using computed tomography (CT-scans) or similar technologies. CT-scans are voxel models and directly describe the volume of an investigated artifact. These can be used as a direct input for a large variety of methods for computational analysis~\cite{Zhang2016a,deGeus2017,Korshunova2020,Korshunova2021a,Gravenkamp2020}.

However, many structures cannot be scanned due to their size or the high attenuation of their material. In such cases, one can resort to optical scanning devices to obtain the shape and topology of existing structures. These scanners deliver oriented point clouds. Oriented point clouds consist of coordinates each of which represents a point on the surface of the scanned structure plus an associated normal vector.  The usual step to obtain an analysis suitable model is then to convert the point cloud to a surface mesh. This surface mesh describes a valid closed volume in the sense of a B-rep model\footnote{i.e. a 2-manifold without a boundary}. The B-rep model is then the starting point of a boundary conforming, volumetric mesh generation. The generated volumetric mesh, which usually consists of (boundary-) conforming tet- or hexahedrals, can then be used as an input to, for example, the finite element method. In many cases it is of course possible to carry out the required pre-processing steps from the acquired shape in form of a point cloud to an analysis suitable mesh. Yet, the generation of such analysis suitable computational models can be expensive and the pipeline is not free from numerous pitfalls. Some of the most common complications are flaws in the surface topology, local over-refinements or other problems due to the tight restrictions on aspect ratios and allowed angles within the mesh, see~\cite{Cottrell2009,Wassermann2019}. 

Recently, it was demonstrated in~\cite{Kudela2020} that is possible to use point clouds as a direct input for computational mechanical analysis and, thereby, to avoid the multitude of otherwise necessary pre-processing steps. Therein, the key idea is to view the point cloud as a geometrical model which implicitly defines the physical domain of interest. We, therefore, discuss the basic aspects of computational mechanics on implicitly defined domains next.

\subsection{Computation on implicitly defined domains}

The computation on implicitly defined geometries is ideally suited for embedded domain methods. Numerous such methods exist since the 1960ies~\cite{Saulev1963}. Some more recent variants include CUTfem~\cite{Burman2015b},  isogeometric B-Rep analysis~\cite{Breitenberger2015}, immersogeometric analysis~\cite{Kamensky2015}, the aggregated unfitted finite element method~\cite{Badia2018}, and the finite cell method (FCM) which is used in the paper at hand~\cite{Duster2008}. We refer to~\cite{Duster2017} for a detailed review including a collection of relevant literature. Interestingly, the FCM has recently also been picked and further developed by other communities as well (see for example contributions in the area of computer graphics~\cite{Longva2020}). 

\begin{figure}[!htbp]
  \includegraphics[width=\textwidth]{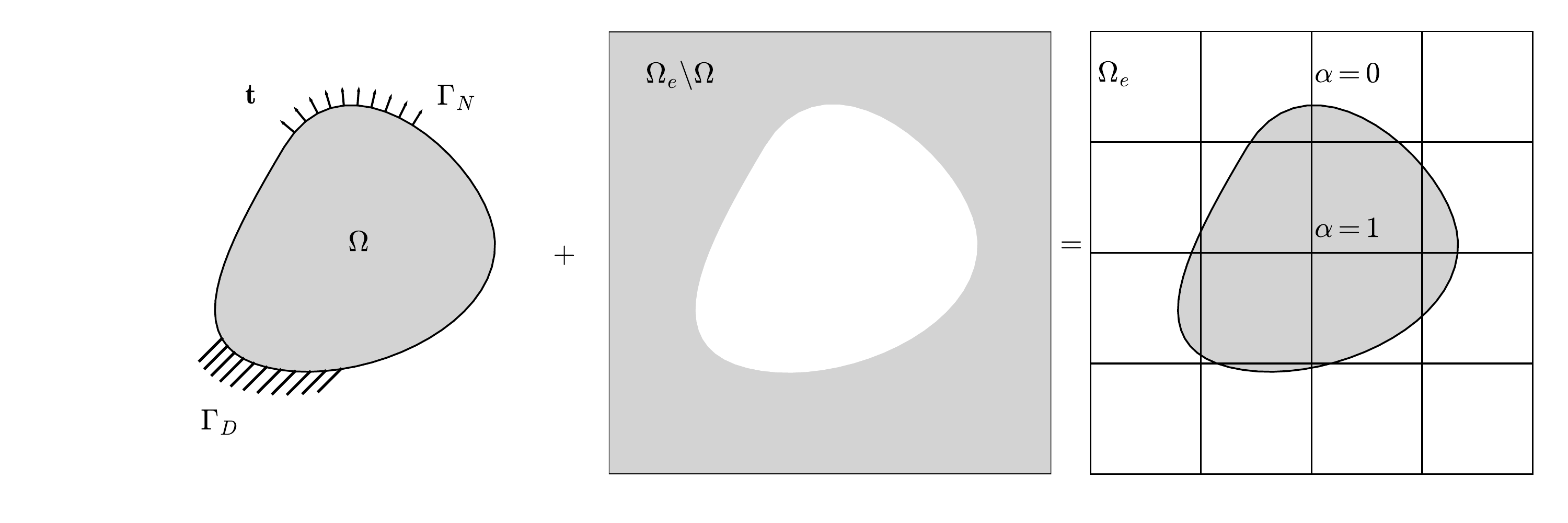}
\caption{The core concept of the FCM. The physical domain $\Omega$ is extended by the fictitious domain $\Omega_e \backslash \Omega$. The extended domain $\Omega_e$ can then be meshed easily. To recover the physical domain $\Omega$, the extended domain $\Omega_e$ is penalized by $\alpha=0$ in $\Omega_e \backslash \Omega$.}
	\label{fig:fcm}
\end{figure}

All of the embedded domain methods mentioned above have different advantages and disadvantages depending on their field of applications but share a common goal: to avoid boundary conforming mesh generation. In the following we present the point of view of the FCM. The FCM extends the physical domain of interest $\Omega$ with a fictitious domain, such that their union $\Omega_e$ can be meshed easily. \Cref{fig:fcm} depicts this basic idea.

To solve the linear elastostatic problem the following weak formulation is then used: ``Find
$\boldsymbol{u} \in H^{1}(\Omega_{e})$ such that

\begin{equation}\label{eq:varPrinc}
  \left(\alpha \, \boldsymbol{\sigma}(\boldsymbol{u}), \nabla \boldsymbol{w}\right)_{\Omega_{e}}
  + \beta \cdot \left(\boldsymbol{u}, \boldsymbol{w} \right)_{\Gamma_{D}}
  = \left(\alpha \, \boldsymbol{b}, \boldsymbol{w} \right)_{\Omega_{e}}
  + \left(\boldsymbol{\hat{t}}, \boldsymbol{w} \right)_{\Gamma_{N}}
  + \beta \cdot \left(\boldsymbol{\hat{u}}, \boldsymbol{w} \right)_{\Gamma_{D}}
  \quad \forall \, \boldsymbol{w} \in H^{1}(\Omega) \, .''
\end{equation}

Here $(\cdot, \cdot),\boldsymbol{\sigma}, \boldsymbol{u}, \boldsymbol{w} \text{ and }
\boldsymbol{b}$ denote the $L^{2}$ scalar product, the Cauchy stress tensor, the
displacement vector, the test function, and the body forces, respectively. For consistency with the original problem formulated on the physical domain, the following indicator function $\alpha$ is defined:
\begin{equation}
  \alpha(\boldsymbol{x})= \left\{\begin{array}{ll}
  {1} & {\forall \boldsymbol{x} \in \Omega} \\
  {\alpha_{fic}} & {\forall \boldsymbol{x} \in
  \Omega_{e} \backslash \Omega}\end{array}\right.
\end{equation}
where $\alpha_{f i c} \ll 1$. Thus, information about the physical domain is accounted
for by the function $\alpha(\boldsymbol{x})$ rather than the mesh itself. To realize
Neumann boundary conditions, in our case the traction $\hat{\boldsymbol{t}}$ on
$\Gamma_{N}$, a contour integral needs to be evaluated. The weak formulation
presented in~\cref{eq:varPrinc} contains two additional contour integrals and a
scalar parameter $\beta$. These stem from the penalty method and enforce the
Dirichlet boundary conditions in the weak sense. Thus, the degree to which
a solution of~\cref{eq:varPrinc} will satisfy the prescribed displacement
$\hat{\boldsymbol{u}}$ on $\Gamma_{D}$ depends on the choice of $\beta$ whereby larger values of
$\beta$ yield a stronger enforcement~\cite{babusFiniteElement}. Other formulations, such as non-symmetric Nitsche versions as presented in~\cite{Schillinger2016a} or discontinuous Galerkin formulations as presented in~\cite{Kollmannsberger2015} would as well be possible. 

In principle, the weak formulation given in~\cref{eq:varPrinc} can be directly applied to various kinds of
geometric models as long as a) the indicator function can be defined such that the volumetric terms can be evaluated accurately and b) as long as the contour integrals can be evaluated on those sections of the boundary where either non-zero Neumann or general Dirichlet boundary conditions are imposed. 

Point cloud models fall into a similar category as they can be generated from pictures or are directly provided by laser scans of existing artifacts. The direct mechanical analysis of objects described by (oriented) point clouds as presented in~\cite{Kudela2020} used a diffuse approach to evaluate the contour integrals in~\cref{eq:varPrinc}. However, a diffuse approach carries mainly two disadvantages: a) it causes a large numerical integration effort close to the diffuse boundary, b) it only provides a clear separation between controlling the solution itself or its gradient on the diffuse boundary in its thin limit. However, the thin limit renders the numerical integration unaffordable thus forcing a trade-off between integration effort and accuracy. 
\vspace{10pt}

\emph{The main novelty of the article at hand is the development and numerical investigation of a sharp interface approach which possesses none of these drawbacks.}

\vspace{10pt}

The structure of the article is as follows:~\Cref{sec:theory} introduces the
theoretical aspects of the diffuse and the new sharp interface
method.~\Cref{sec:numExp} is dedicated to demonstrate the effectiveness of the sharp
interface approach by its application to a benchmark case. Therein, only the more
challenging case of the application of Dirichlet boundary conditions is discussed
before conclusions are drawn in~\cref{sec:conclusions}.


\section{Enforcing essential boundary conditions on domains defined by point clouds} \label{sec:theory}

This section starts by summarizing the challenge of approximating the boundary
locally in~\cref{sec:bound_approx} to clarify the view onto the subject of imposing
boundary conditions in the context of point cloud models and the finite cell method.
The necessity to locally approximate the boundary is common to all presented
approaches. The diffuse interface approach is then rehearsed in~\cref{sec:diffInt}
before the new sharp interface approach is presented in~\cref{sec:sharpInt}. This
section is then concluded by a discussion of principle challenges and limitations
in~\cref{sec:robLim}.

\subsection{Local boundary approximation}\label{sec:bound_approx}

Let $P=\left\{ \boldsymbol{p}_{i} \right\}_{i \in[1, n_{points}]}$ be a point cloud
consisting of points which were sampled from a physical boundary $\Gamma$ with small sampling
errors and equal sampling density. Then the point cloud only holds incomplete and  
noisy information about that boundary. However, each point indicates that the boundary is
located somewhere in its close proximity.

\begin{figure}[!htbp]
\centering
\begin{subfigure}{.29\textwidth}
						\includegraphics[scale=0.5]{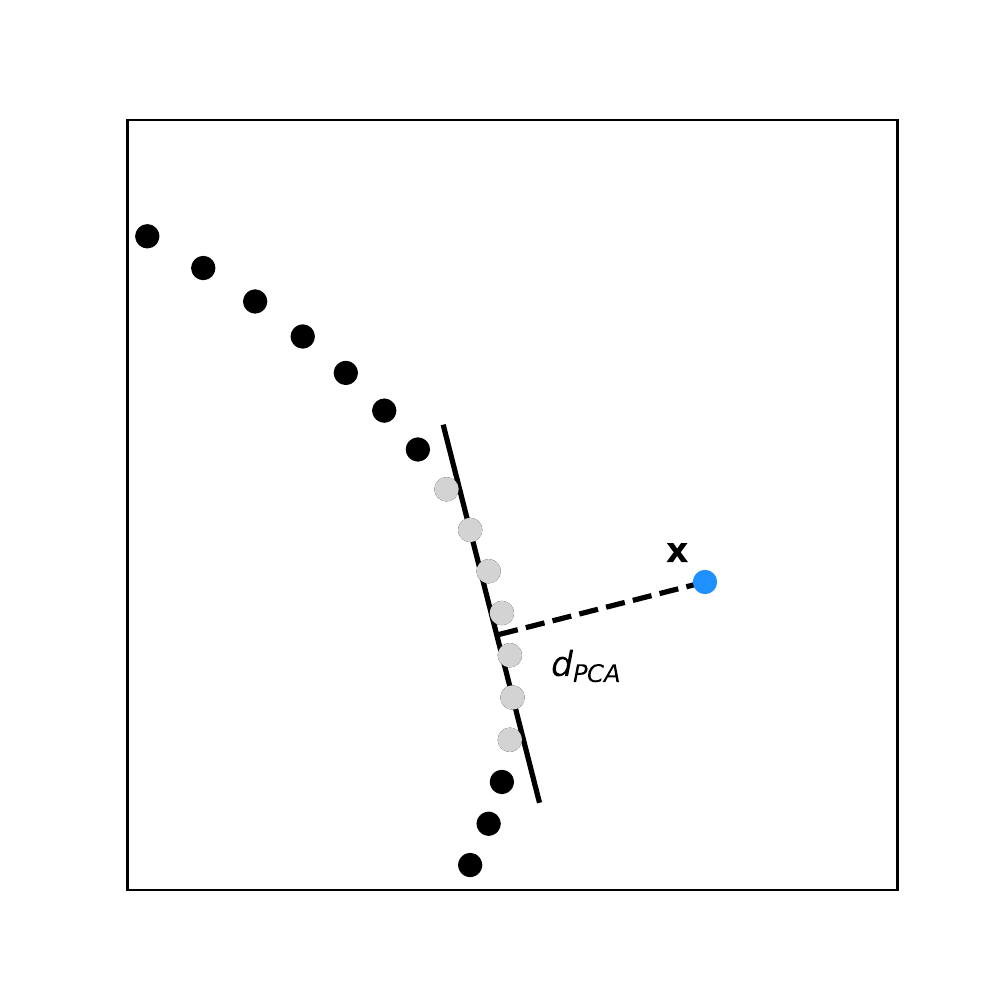}
						\caption{Visualization}
						\label{fig:dpca}
\end{subfigure}						
\begin{subfigure}{.65\textwidth}
\small
		\begin{algorithmic}[1]
	      \Procedure{$d_\text{PCA}$}{$P, \boldsymbol{x}$}
      \State $\{ \boldsymbol{p}_{i}^{k} \} \gets N^{k}(P, \boldsymbol{x})$
      \Comment{Sorted $k$-nearest neighbors}
      \If{$\| \boldsymbol{p}_{1}^{k} - \boldsymbol{x} \| > \boldsymbol{r}$}
      \State \textbf{return} $\|\boldsymbol{p}_{1}^{k}-\boldsymbol{x}\|$
      \Comment{Distance to nearest neighbor}
      \Else
      \State $\overline{ \boldsymbol{p}} \gets \left| \{ \boldsymbol{p}_{i}^{k} \}
      \right|^{-1} \sum_{i}^{k} \boldsymbol{p}_{i}^{k}$
      \Comment{Tangent plane support point}
      \State $ \boldsymbol{C} \gets \sum_{i}^{k}
      \left( \boldsymbol{p}_{i}^{k}-\overline{ \boldsymbol{p}}  \right) \otimes
      \left( \boldsymbol{p}_{i}^{k}-\overline{ \boldsymbol{p}}  \right) \, $
      \Comment{Covariance matrix}
      \State $ \overline{\boldsymbol{n}}_{1}, \overline{\boldsymbol{n}}_{2}
      \gets \text{eigenvectors}(C)$
      \Comment{Tangent plane normal vector}
      \State \textbf{return} $\left| \overline{\boldsymbol{n}}_{2} \cdot
        (\boldsymbol{x} - \overline{\boldsymbol{p}}) \right| $ 
        \Comment{Distance to plane}
      \EndIf
      \EndProcedure
\end{algorithmic}
\caption{Procedure to compute $d_\text{PCA}(\boldsymbol{x})$}
\label{alg:dpca}		
\end{subfigure}			
\caption{ The value assigned to $d_\text{PCA}(\boldsymbol{x})$ is either the test
  point's distance to the least squares fitting plane of its $k$-nearest neighbors $\{
  \boldsymbol{p}_{i}^{k} \}$, $i=1\text{...}k$, or the distance to its nearest
  neighbor $\boldsymbol{p}_{1}^{k}$. The latter is returned if the test point is not
  contained within a circle of radius $\boldsymbol{r}$ around
  $\boldsymbol{p}_{1}^{k}$.}
  \label{fig:dpca_top}
\end{figure}

To construct an implicit approximation of the boundary we assume that for test points
$\boldsymbol{x}$, which are in the vicinity of the unknown boundary $\Gamma$, a local
approximation plane $\Gamma_{\text{PCA}, N^{k}}$ can be computed from the $k$-nearest
neighbors $N^{k}(P, \boldsymbol{x})$~\cite{Kudela2020}.~\Cref{fig:dpca_top}
depicts how a function $d_\text{PCA}(P, \boldsymbol{x})$ giving the unsigned distance
to the approximated boundary can be defined based on this assumption.

Querying $N^{k}(P, \boldsymbol{x})$ can yield results for which
$\Gamma_{\text{PCA},N^k}$ is invalid. This is for example the case if a test point
close to the medial axis of $\Gamma$ is chosen. Points on the medial axis have
multiple sets of closest points on $\Gamma$ and hence the points in $N^{k}(P,
\boldsymbol{x})$ are not necessarily neighbors
themselves~\cite{amentNewVoronoibased}. But this was implied in the construction of
$\Gamma_{\text{PCA},N^k}$ based on $N^{k}(P, \boldsymbol{x})$ as depicted
in~\cref{fig:dpca}. Therefore the parameter $\boldsymbol{r}$ is introduced to limit
the set of points for which $\Gamma_{\text{PCA},N^k}$ is computed.

\subsection{Diffuse interface approach}\label{sec:diffInt}

Assuming $d_{\Gamma_{D}}(\boldsymbol{x})$ is a function giving the unsigned distance
to the Dirichlet boundary then the zero-set of this function is equivalent to the
Dirichlet boundary or $\Gamma_{D}=\{\boldsymbol{x} \, | \,
d_{\Gamma_{D}}(\boldsymbol{x})= 0\}$. In Diffuse Boundary Methods the boundary is
modeled as a layer $\Gamma_{D}^{\epsilon} = \{\boldsymbol{x} \, | \,
d_{\Gamma_{D}}(\boldsymbol{x}) \leq \epsilon\}$ instead, relaxing the sharp interface
condition, with $\epsilon \ll 1$ being a parameter controlling the thickness. The
approach was proposed in the context of approximating partial differential equations
with phase field methods, early works include~\cite{liSOLVINGPDES}
\cite{ratzPDESurfaces}.

Leveraging this model and following our intent to compute contour integrals resulting
from the penalty method on point cloud data we replace the integral over $\Gamma_{D}$
on the left hand side of~\cref{eq:varPrinc} with

\begin{equation}\label{eq:diff_approach}
  \beta \cdot \int_{\Gamma_{D}} \boldsymbol{u} \cdot \boldsymbol{w} \,d \Gamma
  = \beta \cdot \int_{\Omega_{e}} \delta \left( d_{\Gamma_{D}} \right) \left( \boldsymbol{u}
  \cdot \boldsymbol{w} \right) \,d \Omega 
  \approx \beta \cdot \int_{\Omega_e} \delta^{\epsilon} (d_{P C A}(P_{D},
  \boldsymbol{x})) \left( \boldsymbol{u} \cdot \boldsymbol{w} \right)
  \,d \Omega \, .
\end{equation}
\noindent
where
\begin{equation}
  \delta^{\epsilon}(x)=\left\{\begin{array}{ll}{
  \frac{1}{2 \epsilon}\left(1+\cos \left(\frac{\pi x}{\epsilon}\right)\right)}
  & {\text { if }|x| \leq \epsilon} \\ {0} & {\text { otherwise }}\end{array}\right.
\end{equation}
\noindent
Note that this converts the boundary integral into a domain integral by employing the
distance function $d_{\Gamma_{D}}$ and the Dirac delta distribution $\delta(x)$. In
the context of numerical analyses the Dirac delta distribution has to be replaced by
a regularized form $\delta^{\epsilon} (x)$ to make computations feasible, see
e.g.~\cite{leeRegularizedDirac} for details on this replacement. Lastly, since
$d_{\Gamma_{D}}(\boldsymbol{x})$ is not known explicitly but given only by a point
cloud, $d_{\Gamma_{D}}(\boldsymbol{x})$ is replaced with the distance function
$d_\text{PCA}(P_{D}, \boldsymbol{x})$. Therein, $P_{D} \subseteq P$ and $P_{D}$ are
the points belonging to $\Gamma_{D}$. The corresponding boundary integral on the
right hand side of~\cref{eq:varPrinc} is treated accordingly.


Further considerations are necessary for an accurate numerical integration of these
terms. Note that the integrand of the resulting domain integral
in~\cref{eq:diff_approach} is only non-zero in a small fraction of $\Omega_{e}$.
Additionally, the polynomial degree of the integrand is increased due to the
multiplication with $d_\text{PCA}(P_{D}, \boldsymbol{x})$.

\begin{figure}
  \includegraphics[width=\textwidth]{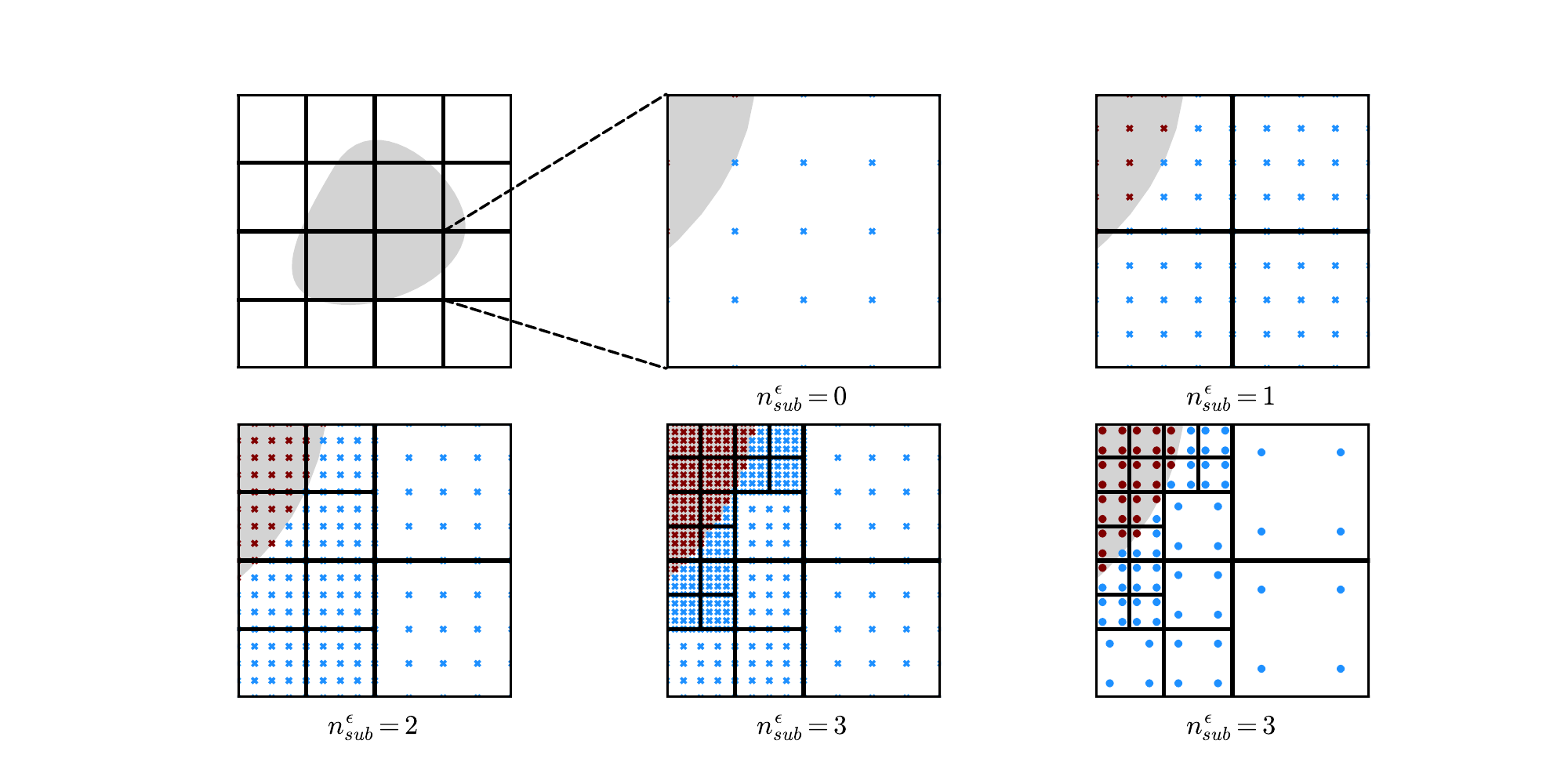}
\caption{Quadtree for diffuse interface approach: 5 equidistant points $\boldsymbol{x}_{i}$ in each direction to test for subdivision until specific subdivision level $n_{sub}^{\epsilon}=3$ is reached. Crosses indicate the test points $\boldsymbol{x}_{i}$, dots the integration points (lowest figure to the right).}
  \label{fig:spacetree}
\end{figure}

Therefore, a straightforward application of standard quadrature rules to integrate
the domain integral containing the diffuse boundary will yield erroneous results. As
a remedy, spacetree-based integration schemes may be applied. These subdivide each
mesh cell recursively until a specified maximum depth $n_{sub}^{\epsilon}$ is
reached. A stopping criterion is defined as follows: subdivide a cell if
$\delta(d_\text{PCA}(P_{D}, \boldsymbol{x}_{i})) > \epsilon_\text{d}$ for any test
point $\boldsymbol{x}_{i}$ from a suitably dense regular grid\footnote{we chose
  $\epsilon_\text{d}=10^{-5}$}. The process is depicted in~\cref{fig:spacetree}.
Subsequently a Gauss-Legendre quadrature of higher order is applied to all subcells
of the deepest recursion depth. The same integration scheme is regularly employed in
FCM to reduce integration errors resulting from the discontinuity of integrands due
to multiplication with $\alpha(\boldsymbol{x})$.

Now we assume the recursive subdivision rule has been applied to a mesh cell giving
$n_{sc}$ subcells. Then the penalty matrix of the cell $\boldsymbol{K}_{P}^\text{ce}$
resulting from discretizing the diffuse term in~\cref{eq:diff_approach} is
approximated as follows

\begin{equation}
  \boldsymbol{K}_{P}^\text{ce} \approx \beta \cdot \sum_{s}^{n_\text{sc}} \sum_{j}^{n_\text{gauss}^{\epsilon}}
  \sum_{i}^{n_\text{gauss}^{\epsilon}}
  \delta^{\epsilon} (d_{P C A}(P_{D}, \boldsymbol{x}(\xi_{i}, \eta_{j})))
  (\boldsymbol{N}^{T} \boldsymbol{N}) \| \boldsymbol{J}^{s} \|
  \| \boldsymbol{J}^{ce} \| |_{\xi_{i}, \eta_{j}} v_{i} v_{j} \, .
\end{equation}

Here $n_\text{gauss}^{\epsilon}, \boldsymbol{N}, \| \boldsymbol{J}^{s} \| \text{ and
} \| \boldsymbol{J}^\text{ce} \| $ denote the number of Gauss-Legendre points in one
dimension, the matrix of shape functions, the Jacobian determinant of the subcell,
and the Jacobian determinant of the cell, respectively. The Gauss-Legendre
coordinates and weights in each dimension are labeled as $\xi_{i}, \eta_{j}$ and
$v_{i}, v_{j}$. The penalty vector $\boldsymbol{f}_{P}^\text{ce}$ on the right hand
side of~\cref{eq:varPrinc} is treated analogously.

The degree to which $\boldsymbol{K}_{P}^\text{ce} \text{ and }
\boldsymbol{f}_{P}^\text{ce}$ enforce the prescribed displacement in the respective
cell depends on a variety of factors. Comparatively large values need to be chosen
for $\beta$ because the penalty method, unlike Nitsche's method, is an inconsistent
method. 


However, a large value of $\beta$ combined with a non-zero $\epsilon$ introduces an unwanted coupling between Dirichlet and Neumann boundary conditions. This additional error must be balanced out by a sufficiently low value of $\beta$ which, in turn, leads to an insufficient approximation of the Dirichlet boundary condition. This dilemma will be clearly pointed out in~\cref{sec:benchmarkAnnularPlate}. Therefore, the parameter $\epsilon$ should be chosen as small as possible since $\delta^{\epsilon}(x)$ converges to $\delta(x)$ only for $\epsilon \rightarrow 0$.  Unfortunately, the lower bound of $\epsilon$ is determined by the effort one wants to spend in the numerical integration. Choosing a smaller $\epsilon$ necessitates to increase the subdivision depth of the spacetree to reach the same accuracy. Increasing $n_\text{sub}^{\epsilon}$ results in exponentially more integration points potentially rendering the numerical integration impractical.  

\subsection{Sharp interface approach}\label{sec:sharpInt}

\begin{figure}[!htbp]
  \centering
   \input{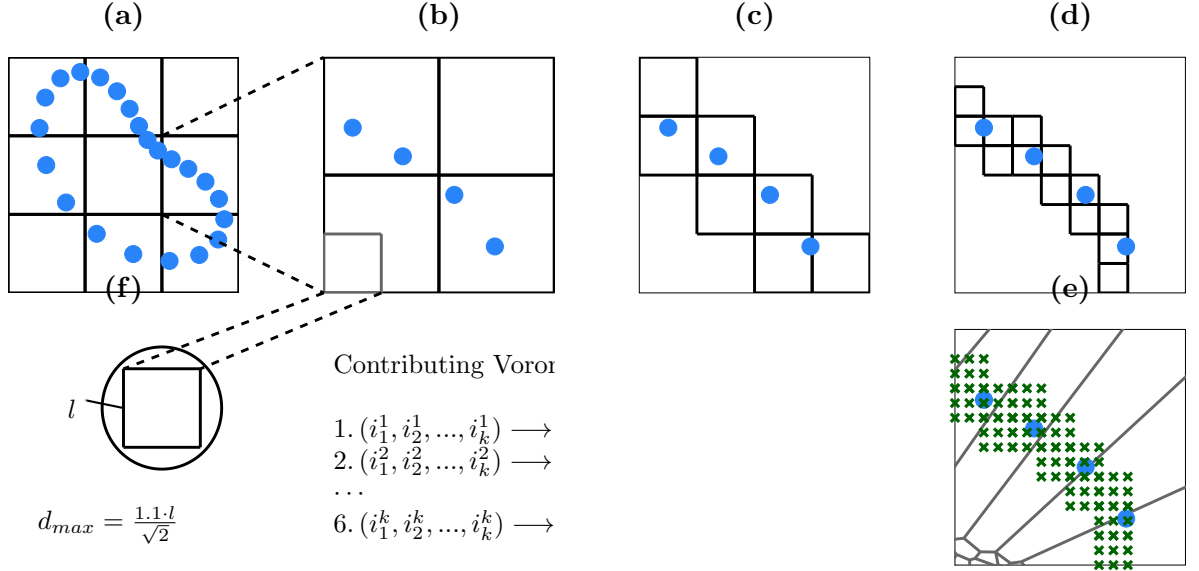}
  \caption{\revb{Step 1: finding contributing Voronoi regions. Each of these regions defines one approximation plane $\Gamma_{\text{PCA},N^{k}}$ approximating one section of the Dirichlet boundary $\Gamma_{D,\text{PCA}}$.}}
  \label{fig:sharpInt1}
\end{figure}

\begin{figure}[!htbp]
  \centering
  \input{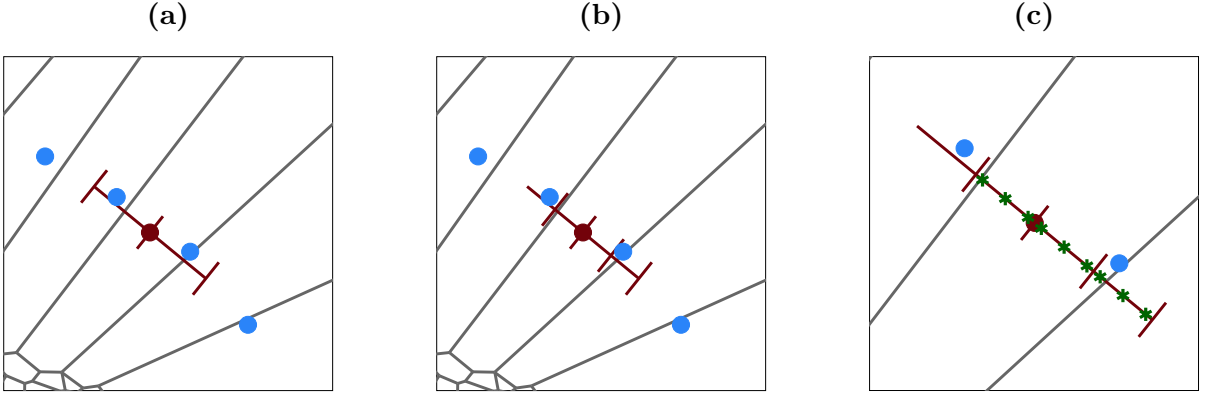}
  \caption{\revb{Step 2: Determining the length of the approximated boundary sections $\Gamma_{\text{PCA},N^{k}}$. Their length is defined by their intersection with the boundaries of the Voronoi regions. Since these boundaries are only implicitly known via querying the distance function $d_\text{PCA}(P_{D},\boldsymbol{x})$ point wise, the intersection (and therefore the length of the boundary sections) is approximated using a bisection algorithm.}}
  \label{fig:sharpInt2}
\end{figure}

\revb{In this section, we propose a different approach which eliminates the problems
  mentioned above by omitting the diffuse boundary assumption. Instead, we propose
  to replace the integral over $\Gamma_{D}$ on the left hand side
  of~\cref{eq:varPrinc} directly with

\begin{equation}
  \beta \cdot \int_{\Gamma_{D}} \boldsymbol{u} \cdot \boldsymbol{w} \,d \Gamma
  \approx \beta \cdot \int_{\Omega_e} \delta (d_{P C A}(P_{D},
  \boldsymbol{x})) \left( \boldsymbol{u} \cdot \boldsymbol{w} \right)
  \,d \Omega 
  = \beta \cdot \int_{\Gamma_{D, P C A}}  \boldsymbol{u} \cdot \boldsymbol{w} 
  \,d \Gamma
  \approx \beta \cdot \int_{\Gamma_{D, sharp}}  \boldsymbol{u} \cdot \boldsymbol{w} 
  \,d \Gamma
	\label{eq:idea}
\end{equation}
\noindent
and treat the corresponding integral on the right hand side of~\cref{eq:varPrinc}
accordingly. Here we propose that it is possible to integrate over
$\Gamma_{D,\text{sharp}}$ a domain which is closely related to the zero-set of
$d_\text{PCA}(P_{D},\boldsymbol{x})$ which we will denote by $\Gamma_{D,\text{PCA}}$.
While this restriction clearly is possible in theory, the crucial question is how an
accurate integration of the right hand side~\cref{eq:varPrinc} can be achieved in a
numerical integration scheme as well. But this direct evaluation seems to be a bold
move only at first sight.

It follows from intuition that a point $\boldsymbol{p}_{i}$ will always be the
nearest neighbor if the test point $\boldsymbol{x}$ is located in its direct
vicinity. In fact every point cloud defines a so-called Voronoi diagram which divides
the respective space into regions where the nearest neighbor point is the same
according to the Euclidean distance~\cite{aurenVoronoiDiagrams-2}. These regions are
commonly denoted as Voronoi regions. A generalization are order-$k$ Voronoi diagrams
which partition the space into order-$k$ Voronoi regions inside which the $k$-nearest
neighbors of each test point $\boldsymbol{x}$ are the same. For any of the subsets of
cardinality $k$ that can possibly be chosen from a given set of points $P$ in
$\mathbf{R}^{2}$ the respective order-$k$ Voronoi region is either a convex, possibly
unbounded polygon or empty~\cite{schmiOrderkVoronoi}. This includes the classical
Voronoi diagram for $k=1$.~\Cref{fig:voronoi} depicts the order-$k$ Voronoi
diagrams of an exemplary point cloud for different $k$.

\begin{figure}[!htbp]
  \centering
  \input{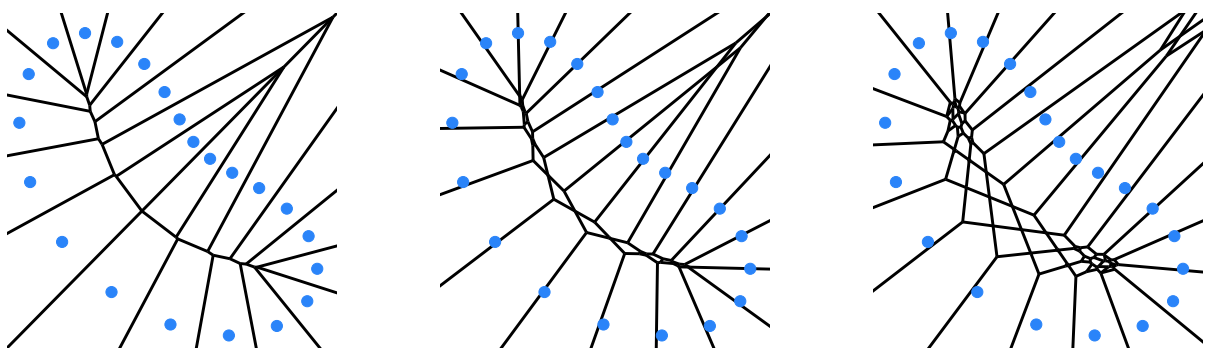}
  \caption{ Order-$k$ Voronoi diagrams of an exemplary point cloud for $k=1,2,4$. }
  \label{fig:voronoi}
\end{figure}

From now on we will refer to Voronoi diagrams and Voronoi regions and imply that they
are generalized of order $k$ for simplicity. The Voronoi diagram provides an
understanding of the reconstructed boundary as given by $d_\text{PCA}(P_{D},
\boldsymbol{x})$. The fact that the $k$-nearest neighbors are the same within each
Voronoi region of the Voronoi diagram implies that the same local approximation plane
$\Gamma_{\text{PCA}, N^{k}}$ is valid in that Voronoi region. Thus, if the Voronoi
diagram of $P_{D}$ is explicitly known, $\Gamma_{D,\text{PCA}}$ could be
reconstructed in $\mathbf{R}^{2}$ as a set of line segments by calculating the points
where each $\Gamma_{\text{PCA},N^{k}}$ intersects with its respective Voronoi region.
Any subsequent numerical integration could, therefore, be carried out directly on
these line segments. Unfortunately this intersection problem is much more complex in
three dimensions for reasons given in~\cref{sec:robLim}. Therefore such an
explicit approach would not translate well to three dimensions.

To this end, we will now introduce a scheme that allows numerical integration over
$\Gamma_{D,\text{sharp}}$ in two dimensions \emph{without having to explicitly compute
Voronoi diagrams}. As such it could be logically extended to three dimensions.
Although $\Gamma_{D,\text{sharp}}$ builds on $\Gamma_{D,\text{PCA}}$ they are in
general not the same. The scheme works on a per mesh cell basis as common in FCM and
consists of two major parts. In the first part depicted in~\cref{fig:sharpInt1} the Voronoi regions that contribute to $\Gamma_{D,\text{PCA}}$ in a respective mesh cell are identified. In the second part depicted in~\cref{fig:sharpInt2} the local approximation plane $\Gamma_{\text{PCA},N^{k}}$ for each contributing
Voronoi region is numerically integrated using a bisection-based scheme.

To identify the Voronoi regions that contribute to $\Gamma_{D,\text{PCA}}$ in a
respective mesh cell the distance function $d_\text{PCA}(P_{D},\boldsymbol{x})$ is
queried in a structured manner. To this end, the mesh cell in question is partitioned
into subcells using a quadtree. A subcell is further subdivided if $d_\text{PCA}
(P_{D}, \boldsymbol{x}_{c}) \leq d_{max}$ holds for the center point of the cell
$\boldsymbol{x}_{c}$. Else the subcell is dismissed since it cannot possibly contain
parts of $\Gamma_{D,\text{PCA}}$. The subdivision process is continued until a
maximum subdivision level $n_{query}^{s}$ is reached. The process to construct this
quadtree is depicted in~\cref{fig:sharpInt1} (a)-(d). The rule to construct $d_{max}$
in a subcell is given in~\cref{fig:sharpInt1} (f).}

\revb{Thus a number of fine grained subcells results for the mesh cell in question. A grid
of test points is distributed on each of these subcells and the $k$-nearest neighbors
of each test point are queried as depicted in~\cref{fig:sharpInt1} (e). The core idea
is that we can understand each query result as answer to the question in which
Voronoi region the test point is located in if the Voronoi regions are identified by
the ordered indices of their $k$ nearest neighbors. This is possible since as we
remember each subset of $k$ points gives one Voronoi region that is either a convex,
possibly unbounded polygon or empty. Hence we can construct a mesh cell specific list
of contributing Voronoi regions by collecting the query results as ordered $k$-tuples
of indices and filtering them for uniqueness. In general the maximum subdivision
level and test grid have to be chosen such that the test points are dense enough to
identify the smallest Voronoi regions that globally contribute to
$\Gamma_{D,\text{PCA}}$.

Since each Voronoi region corresponds to one local approximation plane
$\Gamma_{\text{PCA},N^k}$ each ordered tuple of indices gives a tuple of support
point and normal vector computed as in the distance function
$d_\text{PCA}(P_{D},\boldsymbol{x})$ which is indicated in the center of~\cref{fig:sharpInt1}.
Therefore, the second part of the scheme iterates through the list of intersection planes found in the first part and applies a bisection-based integration scheme to each plane
$(\overline{\boldsymbol{p}}^{j}, \overline{\boldsymbol{n}}_{2}^{j})$. This
bisection-based scheme builds on two central assumptions. Firstly we assume that
$\overline{\boldsymbol{p}}^{j}$ is always inside its respective Voronoi region.
Secondly we assume that a global parameter $l_{max}^{s}$ can be defined to bound the
approximation planes. In 2D this results in each plane
$(\overline{\boldsymbol{p}}^{j}, \overline{\boldsymbol{n}}_{2}^{j})$ being
initialized as a line segment of length $l_{max}^{s}$ as shown in~\cref{fig:sharpInt2} (a). }

\revb{The bisection scheme consists in hierarchically subdividing each initiated line
  segment by querying the $k$-nearest neighbors for a grid of test points. The
  subdivisions continue until a predefined level $n_{sub}^{s}$. Throughout the
  process only subsegments that are either fully within the respective Voronoi region
  or are intersected by the Voronoi regions are kept. Subsegments that are fully
  outside the respective Voronoi regions are neglected since they cannot contain any
  valid integration points (see~\cref{fig:sharpInt2}).

Once these two parts have successfully been executed for a mesh cell giving
$n_\text{vo}$ Voronoi regions each providing a line segment with some
$n^\text{vo,sd}$ subsegments the penalty matrix of that mesh cell
$\boldsymbol{K}_{P}^\text{ce}$ can be approximated as follows

\begin{equation}
  \boldsymbol{K}_{P}^\text{ce} \approx \beta \cdot \sum_\text{vo}^{n_\text{vo}} \sum_\text{sd}^{n_\text{vo,sd}}
  \sum_{i}^{n_{gauss}^{s}} \delta_{\text{ce} \, \text{cell}(\boldsymbol{x}(\xi_{i}))} \,
  \delta_{\text{vo} \, \text{voronoi}(\boldsymbol{x}(\xi_{i}))} \,
  (\boldsymbol{N}^{T} \boldsymbol{N}) \| \boldsymbol{J}^\text{sd} \| |_{\xi_{i}} v_{i} \, .
\end{equation}


Here $n_{gauss}^{s}$ denotes the number of Gauss-Legendre points in one dimension
while $\xi_{i}$ and $v_{i}$ denote the Gauss-Legendre coordinates and weights. The
Kronecker delta function $\delta_{i j}$ is used to check whether an integration point
is contained within the current Voronoi region $\text{vo}$ and current mesh cell
$\text{ce}$. These checks ensure that the numerical integration approximates the
original integral in a consistent manner.}



\subsection{Challenges, limitations and further remarks} \label{sec:robLim}

Both, the diffuse and the sharp interface approach presented in the previous section
build on the fact that the Voronoi regions exhibit regularities for the point clouds
in question. \revb{In 2D for smooth boundaries captured with high and approximately equal
sampling density we can roughly group the Voronoi regions in two groups as visible in
~\cref{fig:voronoiProblem} (a). Firstly there are regions which mainly
extend into normal direction while being confined in tangential
direction~\cite{deyVoronoibasedFeature}. Secondly for $k>1$ 'bands' of small Voronoi
regions form in the vicinity of the medial axis due to fluctuations of $N^{k}(P,
\boldsymbol{x})$. In regions of the first type the edges of the Voronoi regions
should limit the validity of the local approximation plane $\Gamma_{\text{PCA},N^k}$.
Regions of the second type should arguably neither be included in $\Gamma_{D,PCA}$
nor $\Gamma_{D,sharp}$ in. In both approaches $\boldsymbol{r}$ is the control
parameter which can be adjusted to reach these two goals.

For point clouds sampled from non-smooth boundaries and/or with varying sampling
density distinguishing between these two groups of Voronoi regions is no longer
possible. Additionally unfavorable local approximation planes
$\Gamma_{\text{PCA},N^k}$ may exist which leading to discontinuities in
$\Gamma_{D,PCA}$ respectively $\Gamma_{D,sharp}$. This more general setting is
indicated in~\cref{fig:voronoiProblem} (b). This is not to forget that the presence
of noise and outliers could worsen situation significantly and lead to very
unfavorable boundary approximations.

\begin{figure}[!htbp]
  \centering
  \input{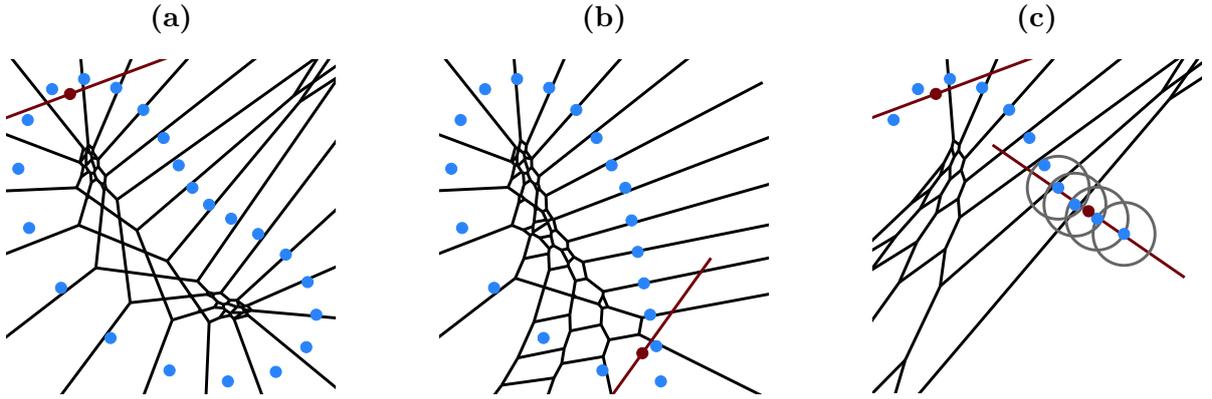}
  \caption{Voronoi diagrams for three different point clouds. Some local
    approximation planes $\Gamma_{\text{PCA},N^k}$ in selected Voronoi regions are
    plotted in red. Circles with favorable radius $\boldsymbol{r}$ displayed
    in gray. }
  \label{fig:voronoiProblem}
\end{figure}

For point clouds sampled from open boundaries approximation planes at the ends could
extend well beyond their intended local reach in the diffuse approach since the no
bound $l_{max}^{s}$ exists. Here it is important to understand that $\boldsymbol{r}$
defines a domain for a point cloud given as the union of disks of radius
$\boldsymbol{r}$ around the individual points $\boldsymbol{p}_{i}$. For small
$\boldsymbol{r}$ this domain truncates $\Gamma_{\text{PCA},N^k}$ in unfavorably
bounded Voronoi regions as indicated in~\cref{fig:voronoiProblem} (c). But
$\boldsymbol{r}$ also has to be large enough such that this truncation does not
introduce unintentional discontinuities in $\Gamma_{D,PCA}$. In a sense the influence
of $\boldsymbol{r}$ in our setting can be compared to the influence of $\alpha$ in
$\alpha$-shapes.

$\alpha$-shapes are related to the union of disks, too, and have been heavily
investigated by the computer graphics community as means to reconstruct a surface
from a point cloud.}~There it became apparent that locally refining $\boldsymbol{r}$
is necessary in many cases ~\cite{edelsShapeSet}. Of course other approaches to the
problem are possible such as introducing dummy points to bound Voronoi regions
~\cite{allieVoronoibasedVariational}. We conclude that it the problem can only
reliably be solved with complementary local parameters, especially in the presence of
noise. However, local parameters would add many degrees of freedom to the overall
approximation $\Gamma_{D,PCA}$. To not complicate the analysis further we opted
against introducing more parameters. Moreover we assert that the central deficiency
of both approaches in their current version is that they cannot in general be applied
to point clouds sampled from open or non-smooth boundaries and restrict the numerical
studies accordingly. However, a more detailed study is beyond the scope of this
paper. Instead we would like to point out that 2D point cloud curve reconstruction is
indeed an active field of research and direct the interested reader to the very
recent state of the art review published in~\cite{Ohrhallinger2021}.


To conclude this section we will briefly discuss the situation in three dimensions.
In three dimensions, Voronoi regions are convex, possibly unbounded polyhedrons.
Efficient explicit computations with the set of points, edges and polygons that make up these
3D Voronoi diagrams would necessitate further considerations in algorithms and
data structures. Additionally, the computation of Voronoi diagrams is challenging
as the diagrams have worst case exponential complexity \cite{aurenVoronoiDiagrams-2}.

However, the presented sharp interface approach does not need the explicit computation of Voronoi diagrams. Instead, the presented approach relies on an implicit representation of the zero-level set of $d_\text{PCA}(P_{D},\boldsymbol{x})$ or $\Gamma_{D,\text{PCA}}$. These would be found by octrees. Quadtrees can then be used to perform the integrations over the boundary planes $\Gamma_{\text{PCA},N^k}$. We would further like to remark that the implicitly reconstructed boundary planes do not even have to possess a $C^0$ continuity because models with flawed boundaries may also be handled by the FCM~\cite{Wassermann2019}. \reva{Although these are all promising hints that an extension to three dimensions might be straight forward, experience shows that only an implementation for 3D will disclose possible pitfalls.}

%
%

%


\section{Numerical examples} \label{sec:numExp} 

\revb{In this section we will investigate two types of examples.~\Cref{sec:benchmarkAnnularPlate} closely investigates a benchmark example known in the community of computational mechanics at which the diffuse interface approach rehearsed in~\cref{sec:diffInt} is compared to the newly developed sharp interface approach presented in~\cref{sec:sharpInt}. All algorithms are implemented in the FCMLab framework which documented in~\cite{Zander2014} and publicly accessible at \url{https://gitlab.lrz.de/cie_sam_public/fcmlab/} 

  In~\cref{sec:PointCloudExamples} we then apply the sharp interface approach on two
  point cloud examples taken from a benchmark-series of point
  clouds~\cite{Ohrhallinger2021} which was built for the evaluation of reconstruction
  algorithms to demonstrate that the algorithm is capable of capturing also examples
  with convex, concave and non-smoothly boundaries involving kinks.}

\subsection{Annular plate benchmarks} \label{sec:benchmarkAnnularPlate}

\subsubsection{Problem definition} \label{sec:sphereWithHole}

\begin{figure}[!htbp]
  \centering
  \includegraphics{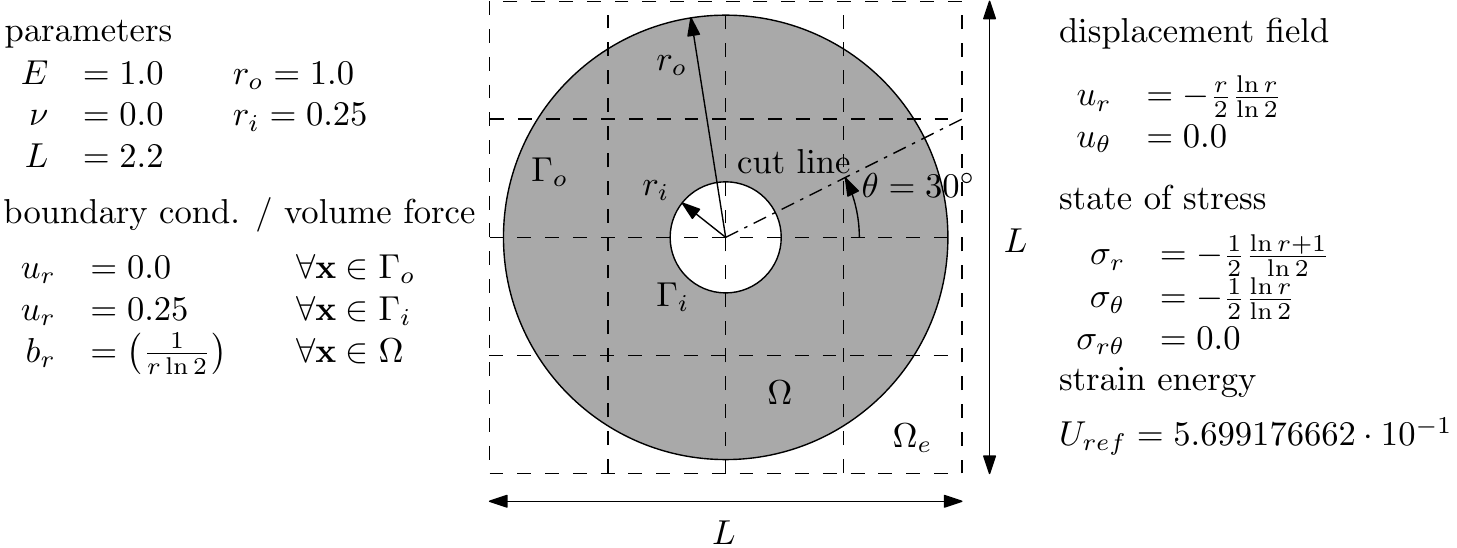}
  \caption{Annular plate setup}
  \label{fig:annular}
\end{figure}

\begin{figure}[!htbp]
  \centering
    \includegraphics[scale=1.0]{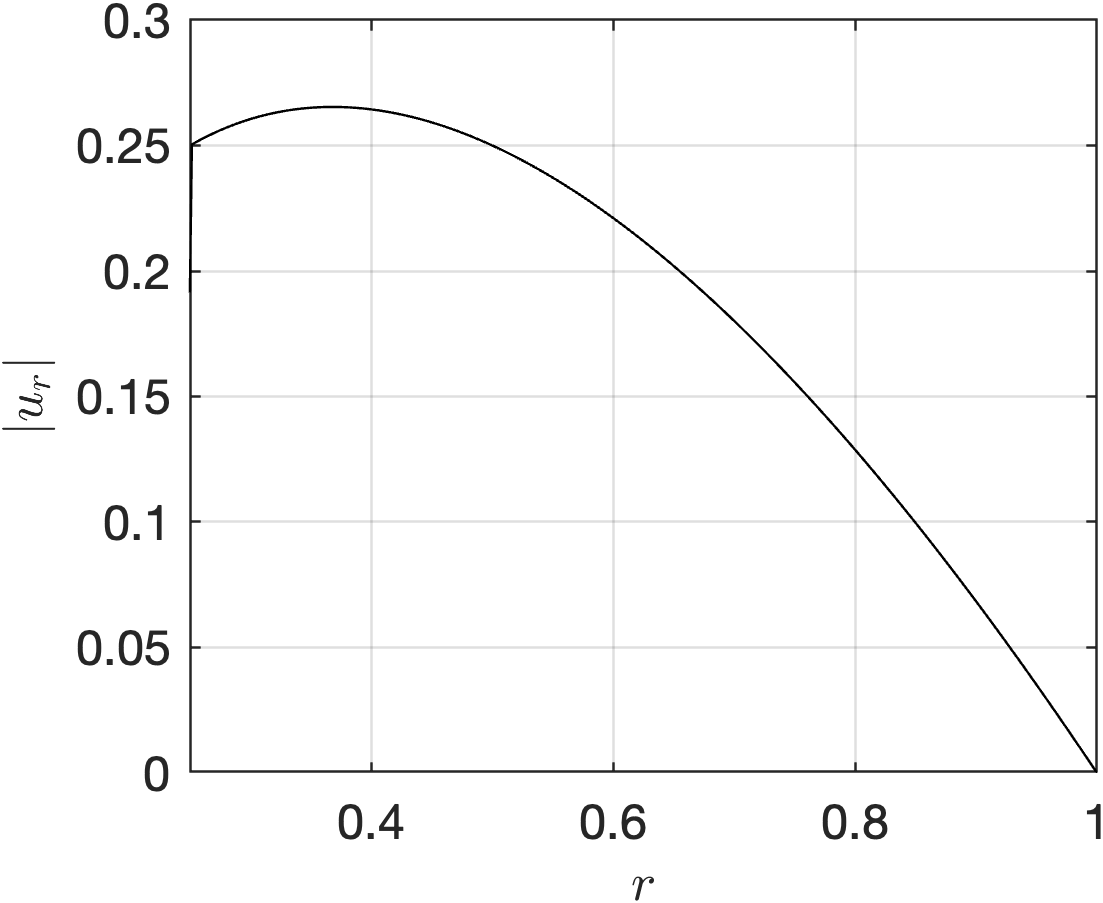}
    \caption{Solution along the radius.}
  \label{fig:solutionUr}
\end{figure}

An annular plate under plane stress conditions with prescribed displacement along the
whole boundary is subjected to a body force in radial direction. A detailed
illustration of the problem, mesh and analytical solution as given in
\cite{zandeFiniteCell} is depicted in~\cref{fig:annular}. A plot of the
solution $u_r$ along the cut line is given in~\cref{fig:solutionUr}. 

Irrespective of the way the penalty integrals in~\cref{eq:varPrinc} are
evaluated the method's effectiveness crucially depends on the choice of $\beta$. From
other works~\cite{nguyeDiffuseNitsche} it is also clear that the diffuse interface approach is very sensitive to the choice of $\epsilon$.
Hence a numerical study is now considered that allows to characterize the
influence of $\beta$ on the resulting error while considering $\epsilon$ in the
diffuse interface approach. The employed error measure is the relative error in the
energy norm defined as

\begin{equation}
  e=\sqrt{\frac{\left|U_{n u m}-U_{r e f}\right|}{U_{r e f}}} \cdot 100 \,[\%] \, .
\end{equation}

The mass matrix and load vector for basis functions of polynomial degree 10 were
precomputed employing a quadtree of depth 10 and a Gaussian quadrature of 11 points.
In this case the indicator function $\alpha(\boldsymbol{x})$, with
$\alpha_{fic}=10^{-8}$, is constructed from exact circle formulas to minimize the
integration error. These precomputations allow to apply the penalty integrals using
the two approaches to the same discrete matrix-vector pair which further enhances their
comparability.\footnote{In line with a reviewer remark we would like to comment that it is of course possible to choose any other type of discretization also combined with local refinements as, e.g. presented in~\cite{Zander2015,DAngella2018,Kopp2021a} but opted not to do so out of simplicity.}

The two point clouds discretizing the boundary with approximately equal sampling
density are constructed by taking $n_{points}=10^{4}$ points on the inner circle and
$4n_{points}$ points on the outer circle. We chose to fix $k=4$ at all times for the
computation of the local approximation planes. For the diffuse interface approach we
restrict $\epsilon$ to the values $\epsilon_{1}=5 \cdot 10^{-3}$, $\epsilon_{2}=5
\cdot 10^{-4}$ and $\epsilon_{3}=5 \cdot 10^{-5}$. 
Further, we set $n_{gauss}^{\epsilon}=10$ and $n_{sub}^{\epsilon}=7$ for $\epsilon_{1}$, respectively
$n_{sub}^{\epsilon}=10$ and $n_{sub}^{\epsilon}=13$ for $\epsilon_{2}$ and
$\epsilon_{3}$. While for the diffuse interface approach the integration parameters
have to be chosen with respect to $\epsilon$ the sampling density is decisive for the
sharp interface approach. Thus we chose $n_{query}^{s}=12$, $n_{sub}^{s}=3$,
$n_{gauss}^{s}=11$ and $l_{max}^{s}=3 \cdot 10^{-4}$ in accordance with the problem's
sampling density.

\revb{To compute a numerical reference the inner circle is discretized with $10^{4}$ line
segments and the outer circle with $4 \cdot 10^{4}$ line segments respectively. The
penalty integrals in~\cref{eq:varPrinc} are then evaluated by computing a
Gaussian quadrature with 11 points on each line segment.}

\subsubsection{Numerical Results} \label{sec:numSecRes}

\begin{figure}[!htbp]
  \centering
%
%
\begin{tikzpicture}

\begin{axis}[%
width=4.844in,
height=3.379in,
at={(0.812in,0.475in)},
scale only axis,
xmode=log,
xmin=10,
xmax=10000000,
xminorticks=true,
xlabel style={font=\color{white!15!black}},
xlabel={Penalty parameter $\beta$},
ymode=log,
ymin=0.01,
ymax=100,
yminorticks=false,
ylabel style={font=\color{white!15!black}},
ylabel={Relative error in energy norm $e$ [\%]},
axis background/.style={fill=white},
xmajorgrids,
xminorgrids,
ymajorgrids,
yminorgrids,
legend style={at={(0.03,0.03)}, anchor=south west, legend cell align=left, align=left, draw=white!15!black, fill=white}
]
\addplot [color=black, line width=1.0pt, mark size=2.3pt, mark=triangle, mark options={solid, black}]
  table[row sep=crcr]{%
10.7721734501594	10.6762804886223\\
17.9690683190231	2.57384329389637\\
29.9742125159471	4.95747309177503\\
50	5.5304022813078\\
83.4050268600029	5.63077365778708\\
139.127970110356	5.9457974421128\\
232.079441680639	6.62912637167555\\
387.131841340563	7.59389063881065\\
645.774832507442	8.66219554561089\\
1077.21734501594	9.65722726859051\\
1796.90683190231	10.4679438060813\\
2997.4212515947	11.0855162517535\\
5000	11.5879442358421\\
8340.50268600029	12.0905417056802\\
13912.7970110356	12.7038966736052\\
23207.9441680639	13.5108205375076\\
38713.1841340564	14.5472942582085\\
64577.4832507441	15.7800758946179\\
107721.734501594	17.0971032808054\\
179690.683190231	18.3339465415085\\
299742.12515947	19.3410255366022\\
500000	20.0560634658808\\
834050.268600029	20.5156873498078\\
1391279.70110356	20.5695837114296\\
2320794.41680639	20.9412683027362\\
3871318.41340564	20.7395289795708\\
};
\addlegendentry{$\epsilon{}_{\text{1}}\text{ = 5.0e-03}$}

\addplot [color=black, line width=1.0pt, mark size=2.5pt, mark=square, mark options={solid, black}]
  table[row sep=crcr]{%
10.7721734501594	10.7438937629196\\
17.9690683190231	3.08632809274439\\
29.9742125159471	4.39614364733924\\
50	4.6458299269065\\
83.4050268600029	4.10465288022073\\
139.127970110356	3.41936527604685\\
232.079441680639	2.79505003571522\\
387.131841340563	2.31047648892095\\
645.774832507442	2.01008612978102\\
1077.21734501594	1.93018547761736\\
1796.90683190231	2.086490633571\\
2997.4212515947	2.46470654452918\\
5000	3.04245076655004\\
8340.50268600029	3.80678316791652\\
13912.7970110356	4.74762556385623\\
23207.9441680639	5.83873031706183\\
38713.1841340564	7.0176960461218\\
64577.4832507441	8.17834012008216\\
107721.734501594	9.19351223076389\\
179690.683190231	9.96715248051698\\
299742.12515947	10.4745847309399\\
500000	10.7335711170405\\
834050.268600029	10.9148619086779\\
1391279.70110356	10.952101022812\\
2320794.41680639	11.005455180536\\
3871318.41340564	10.8891268807944\\
};
\addlegendentry{$\epsilon{}_{\text{2}}\text{ = 5.0e-04}$}

\addplot [color=black, line width=1.0pt, mark size=3.5pt, mark=o, mark options={solid, black}]
  table[row sep=crcr]{%
10.7721734501594	10.7446904695101\\
17.9690683190231	3.0912828124496\\
29.9742125159471	4.39011144320037\\
50	4.63610052532601\\
83.4050268600029	4.086034874304\\
139.127970110356	3.3817068092669\\
232.079441680639	2.71728541994542\\
387.131841340563	2.15004006374632\\
645.774832507442	1.68775525292572\\
1077.21734501594	1.32127095922554\\
1796.90683190231	1.03792158004749\\
2997.4212515947	0.827771520775556\\
5000	0.687157234247418\\
8340.50268600029	0.619610927670521\\
13912.7970110356	0.630620304595913\\
23207.9441680639	0.719214219977899\\
38713.1841340564	0.879852227502146\\
64577.4832507441	1.1113993061909\\
107721.734501594	1.42082591282052\\
179690.683190231	1.78240599389191\\
299742.12515947	2.33050648840392\\
500000	2.97936431656711\\
834050.268600029	3.77597429439525\\
1391279.70110356	4.72904730101228\\
2320794.41680639	5.76805110414651\\
3871318.41340564	6.31933125183495\\
};
\addlegendentry{$\epsilon{}_{\text{3}}\text{ = 5.0e-05}$}

\addplot [color=black, line width=1.0pt, mark size=3.5pt, mark=asterisk, mark options={solid, black}]
  table[row sep=crcr]{%
10.7721734501594	10.7757817590174\\
17.9690683190231	3.12766524782193\\
29.9742125159471	4.38360829161717\\
50	4.63589976095199\\
83.4050268600029	4.08737065047249\\
139.127970110356	3.38316321193073\\
232.079441680639	2.71820077899153\\
387.131841340563	2.1498097949584\\
645.774832507442	1.68534879854152\\
1077.21734501594	1.31453760624121\\
1796.90683190231	1.02226733663657\\
2997.4212515947	0.793585784312208\\
5000	0.615421823295225\\
8340.50268600029	0.476965117864401\\
13912.7970110356	0.369525408277005\\
23207.9441680639	0.286224278535863\\
38713.1841340564	0.221663323667915\\
64577.4832507441	0.171629998723373\\
107721.734501594	0.132847104373528\\
179690.683190231	0.102765022497382\\
299742.12515947	0.0794256800434878\\
500000	0.0612851888368104\\
834050.268600029	0.0470598171552557\\
1391279.70110356	0.0360554554251943\\
2320794.41680639	0.0153113614840775\\
3871318.41340564	0.0187356810513978\\
};
\addlegendentry{Sharp interface}

\addplot [color=black, line width=1.0pt, mark size=3.5pt, mark=x, mark options={solid, black}]
  table[row sep=crcr]{%
10.7721734501594	10.7446907889521\\
17.9690683190231	3.0912834870957\\
29.9742125159471	4.390094906722\\
50	4.63604977877787\\
83.4050268600029	4.08590417101716\\
139.127970110356	3.38139659599556\\
232.079441680639	2.71657950617323\\
387.131841340563	2.14847054830473\\
645.774832507442	1.68431259084187\\
1077.21734501594	1.3137859790704\\
1796.90683190231	1.02177366253707\\
2997.4212515947	0.793328099005429\\
5000	0.615387925858892\\
8340.50268600029	0.477155099189954\\
13912.7970110356	0.369952509551246\\
23207.9441680639	0.28691636747013\\
38713.1841340564	0.222666824022096\\
64577.4832507441	0.173006707393526\\
107721.734501594	0.134687870218227\\
179690.683190231	0.105175665469256\\
299742.12515947	0.0825783864384607\\
500000	0.0653277848668877\\
834050.268600029	0.0523032582942253\\
1391279.70110356	0.0421765144103316\\
2320794.41680639	0.106737610516816\\
3871318.41340564	0.0396681017510703\\
};
\addlegendentry{Reference}

\end{axis}
\end{tikzpicture}%
  \caption{$\beta$-study}
  \label{fig:optBeta}
\end{figure}
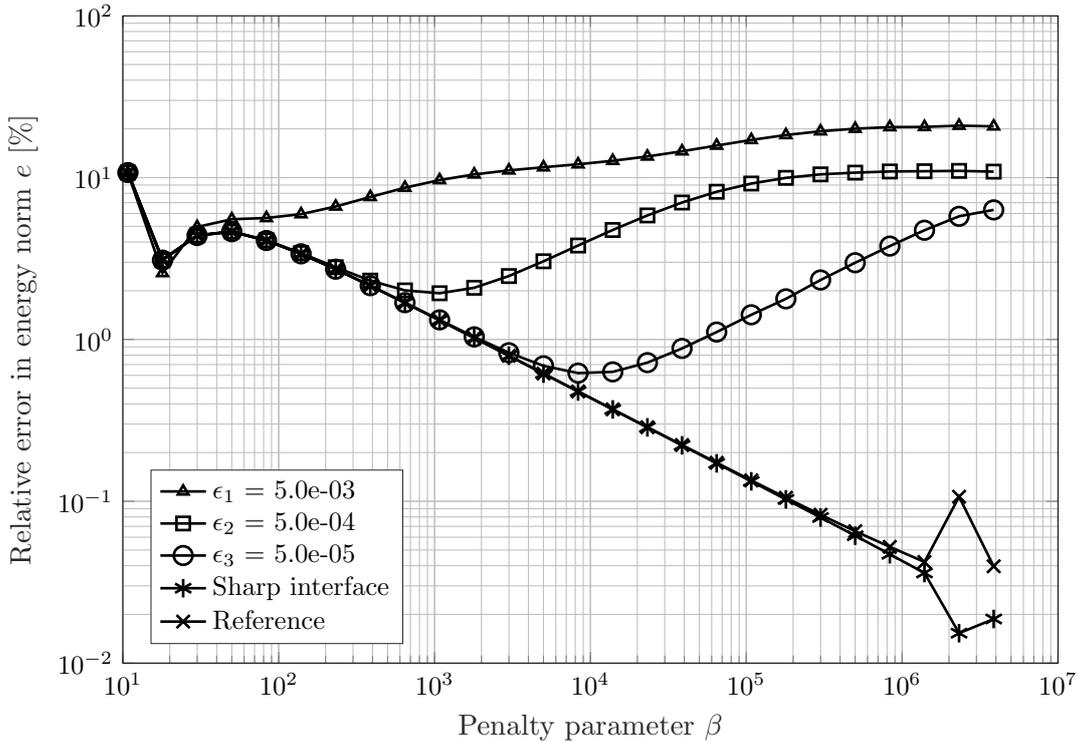

Building on this setup a $\beta$-study is conducted for 26 possible values of $\beta$
logarithmically spaced between 10 and $5 \cdot 10^{6}$. The resulting error $e$ for
the diffuse interface approach, the sharp interface approach and the numerical
reference is plotted in~\cref{fig:optBeta}.

Contrary to the expectation the error $e$ does not monotonically decrease for the diffuse interface approach when larger values $\beta$ are employed. While the error
corresponds to the error of the numerical reference for small $\beta$ the error
curves level off at certain points and the error increases again. Thus the error will
only become minimal for certain optimal values of $\beta$. These optimal values depend on the width of the diffuse boundary $\epsilon$. By contrast, the numerical solution for the sharp
interface approach converges towards the exact solution at the same rate as the
reference solution. The best result of $e \approx 0.04$ is obtained for
the sharp interface approach with $\beta = 1.39 \cdot 10^{6}$, an accuracy which cannot be reached by the diffuse interface approach.

The reason for the convergence behavior of the diffuse interface approach is
examined further next. The point clouds discretizing the inner and outer circle were
chosen with a very high sampling density. Therefore a geometric mismatch between
$\Gamma_{D}$ and $\Gamma_{D,PCA}$ can be ruled out as primary cause of the error.
This is also supported by the results for the sharp interface approach since it
integrates over $\Gamma_{D,sharp}\approx\Gamma_{D,PCA}$.

\begin{figure}[!htbp]
  \centering
  \includegraphics[scale=1.0]{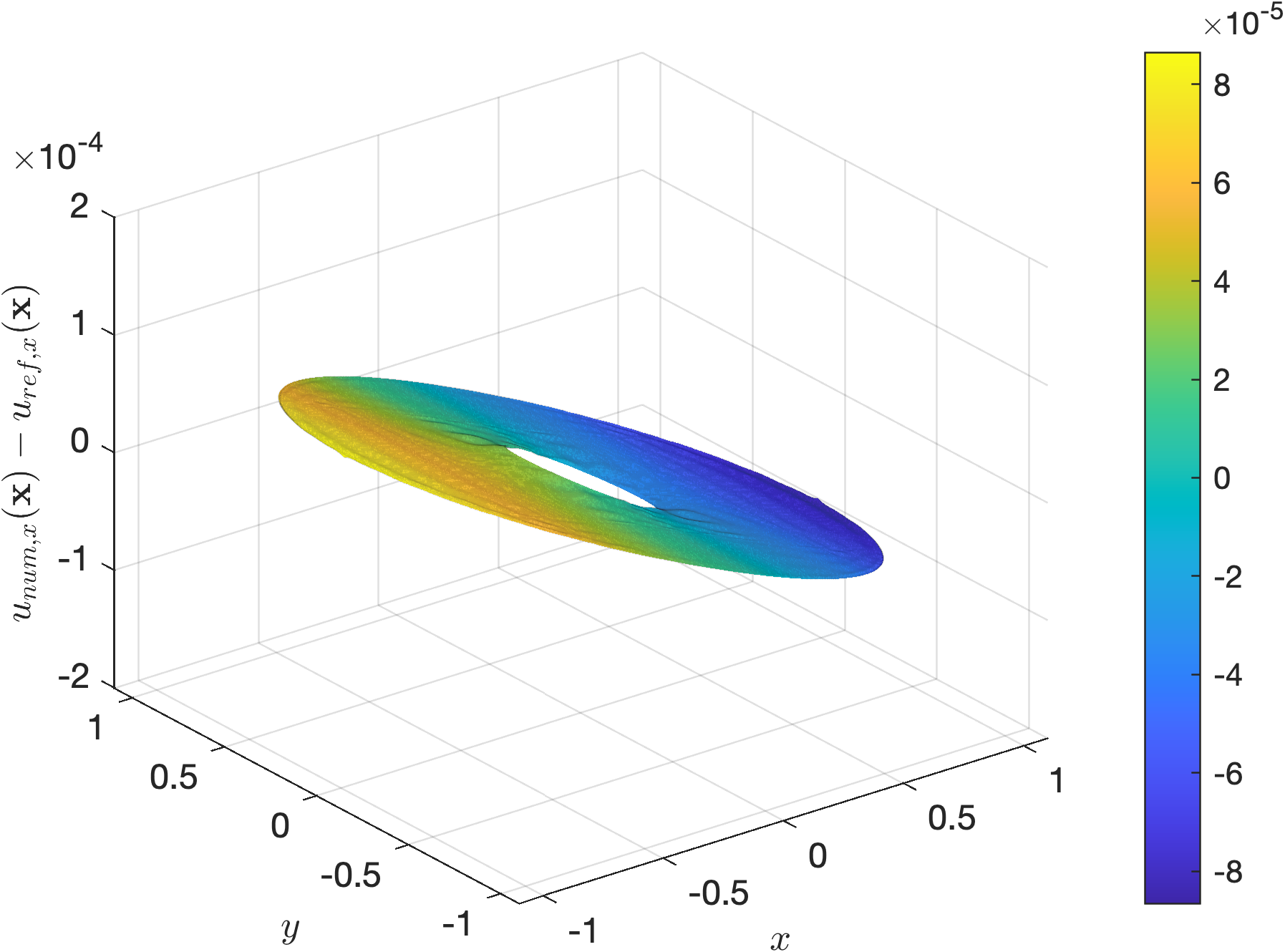}
  \caption{\revb{Abs. error in $u_x$ for the diffuse interface approach with $\epsilon_{3}$
    and $\beta = 8043$}}
  \label{fig:optBetaError1}
\end{figure}
\begin{figure}[!htbp]
  \centering
  \includegraphics[scale=1.0]{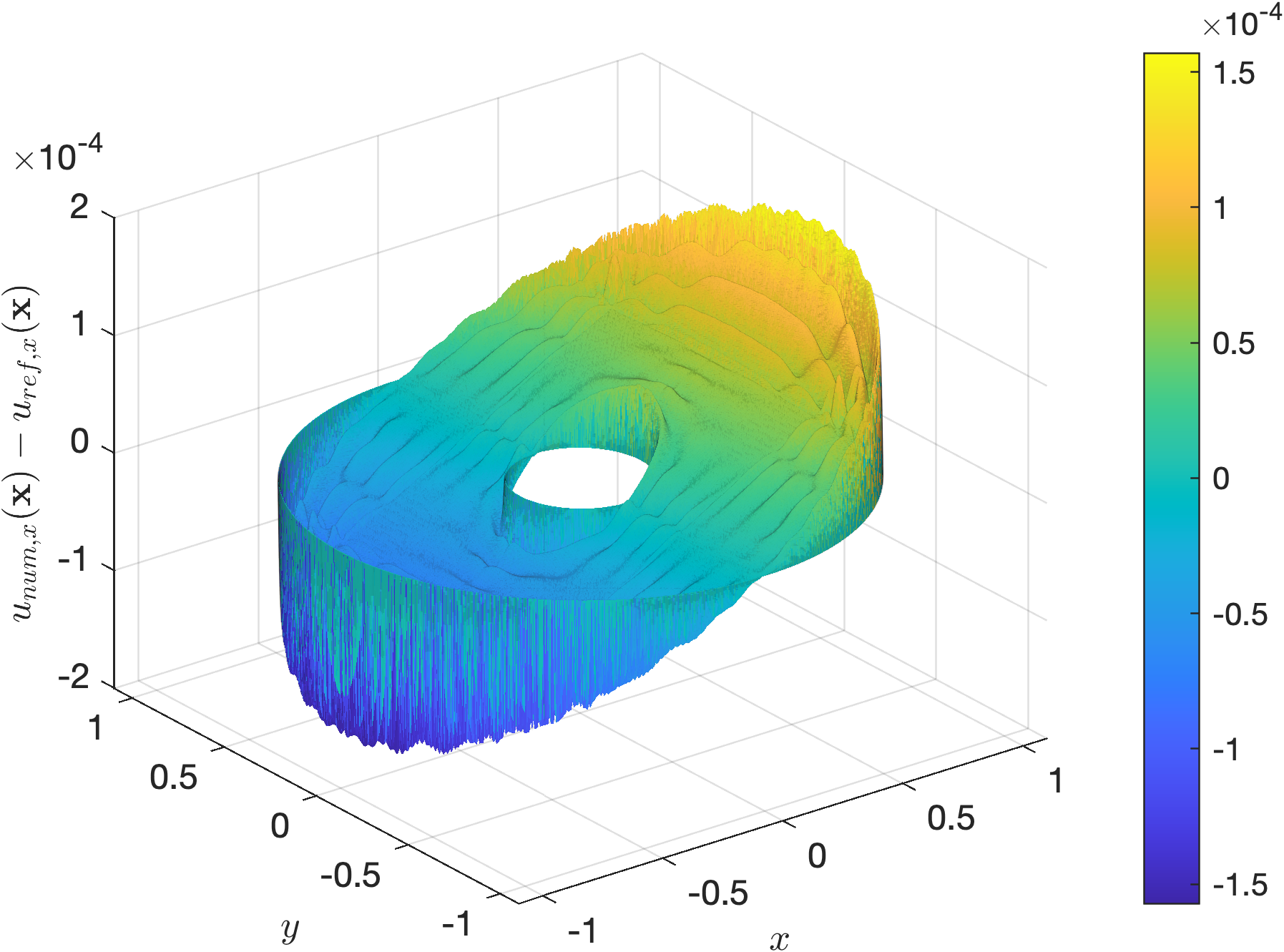}
  \caption{\revb{Abs. error in $u_x$ for the diffuse interface approach with $\epsilon_{3}$
    and $\beta = 5 \cdot 10^{5}$}}
  \label{fig:optBetaError2}
\end{figure}

But no matter how the penalty integrals are computed $\beta$ controls the extent to
which deviations from the constraint are penalized. Therefore, the attainable error
below the tear off point is limited by $\beta$ itself. For large $\beta$ the
constraint is enforced to some extent on the diffuse layer $\Gamma_{D}^{\epsilon}$,
too. In other words, controlling the Dirichlet boundary condition is also constraining the Neumann boundary condition. This undesirable constraining of the gradient towards zero in normal
direction as exemplarily depicted in~\cref{fig:optBetaError2}. This unintended
side effect is stronger for larger $\epsilon$ and larger $\beta$. Due to this
situation a local minimum of the error $e$ in between small and large $\beta$ exists
if $\epsilon$ is small enough. This situation is depicted for $\epsilon_{3}$ in
~\cref{fig:optBetaError1}. The finding that in a diffuse boundary approach only one problem dependent choice of $\beta$ exists for small enough $\epsilon$ is consistent with results by other researchers
\cite{nguyeDiffuseNitsche} \cite{burgeAnalysisDiffuse}.

Very small values of $\epsilon$ are necessary to minimize constraining the gradient. This, in turn, requires costlier numerical integration with a higher
subdivision level $n_{sub}^{\epsilon}$. A trade-off between integration effort and
accuracy results in the fact that the diffuse interface approach is not only less accurate but also computationally more expensive compared to the sharp interface method. Moreover, the range of $\beta$ values leading to favorable results is not known beforehand given a specific problem and a value of
$\epsilon$. 

\begin{figure}[!htbp]
  \centering
  \includegraphics[scale=1.0]{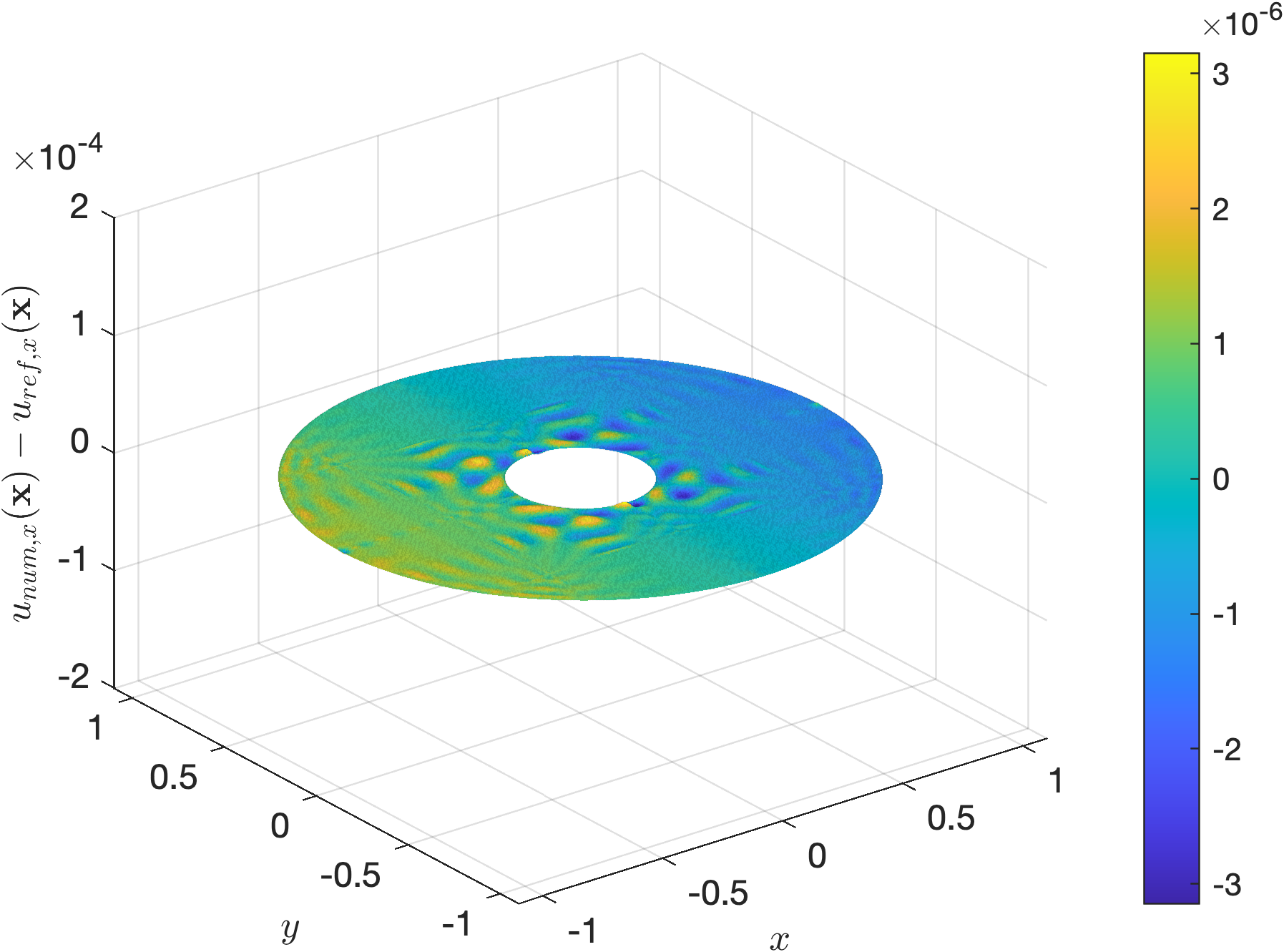}
  \caption{\revb{Abs. error in $u_x$ for the sharp interface approach with $\beta = 5 \cdot
    10^{5}$}}
  \label{fig:optBetaError3}
\end{figure}

The side effect of constraining the gradient does not exist for the sharp interface
approach by construction. Larger values for $\beta$ consistently yield lower errors.
\Cref{fig:optBetaError3} exemplarily depicts the error in $u_x$. While both
approaches start with the same ad hoc boundary approximation $d_{PCA}(P,
\boldsymbol{x})$ they diverge in the way integration points are distributed in the
vicinity of $\Gamma_{D,PCA}$. In this example Gauss quadrature is applied to 2D
sub-cells of quadtrees for the diffuse interface approach. In the sharp interface
approach Gauss quadrature is applied to line segments resulting in much less
integration points. Hence the sharp interface approach does not only evade the trade
off between integration effort and accuracy but implies a lower computational effort
by construction.

\subsection{Point clouds resembling concave, convex and non-smooth boundaries} \label{sec:PointCloudExamples}

\revb{To study the practicability of the sharp interface approach we investigate a
  membrane structure with prescribed deflection on the approximated boundary
  $\Gamma_{D,\text{sharp}}$. The problem can be modeled with Poisson's equation.
  Incorporating the prescribed deflection via the penalty method leads to the
  following weak formulation ``Find $u \in H^{1}(\Omega)$ such that

\begin{equation}
  (\nabla u,\nabla w)_{\Omega} + \beta \cdot (u,w)_{\Gamma_{D,\text{sharp}}}
  = -(\frac{p}{\sigma_0t},w)_{\Omega} + \beta \cdot (u,\hat{u})_{\Gamma_{D,\text{sharp}}}
  \quad \forall \, w \in H^{1}(\Omega) \, .''
\end{equation}
Here $u, w$ and $\hat{u}$ denote the deflection, the test function and
the prescribed deflection, respectively. The vertical load, prestress and thickness are given by
$p, \sigma_0$ and $t$.

The point clouds in question are taken from the benchmark by Ohrhallinger et al.
~\cite{Ohrhallinger2021} and may be found at
\url{https://gitlab.com/stefango74/curve-benchmark}. The point clouds are scaled such
that they fit into the solution domain $[-1.1,1.1]\times[-1.1,1.1]$ which is discretized by
16x16 elements of polynomial degree 10. The solution is fixed to zero along the mesh
boundary. Further we chose $\beta=1\cdot10^{6}$, $\frac{p}{\sigma_0t}=10$,
$\hat{u}=1$, $k=4$ and employ a Gaussian quadrature of 11 points throughout all parts of the example.

The first example \texttt{mc4.txt} consists of multiple convex and concave curves and is depicted in~\cref{fig:curve1a}.
Since all curves are smooth, closed and sampled with approximately equal sampling density
the point cloud is favorable with regard to~\cref{sec:robLim}. Nevertheless the
example poses a challenge due to the structure of the Voronoi diagram. The Voronoi
regions whose approximation planes $\Gamma_{PCA,N^{k}}$ we want to include extend in
normal direction while being confined in tangential direction as visible in
\cref{fig:curve1a}. Due to fluctuations of $N^{k}(P, \boldsymbol{x})$ 'bands' of
small Voronoi regions form in the vicinity of the medial axis which we do not want to
include in $\Gamma_{D,\text{sharp}}$.

Thus the challenge is to select the radius $\boldsymbol{r}$ just right to include all
Voronoi regions in the former group while excluding those in the latter. We find that
this division is in fact possible for $\boldsymbol{r}=0.02$ and a query quadtree of
depth $n_{query}^{s}=5$ per mesh cell. The approximation planes $\Gamma_{PCA,N^{k}}$
initialized with $l_{max}^{s}=8 \cdot 10^{-2}$ and bisected to level $n_{sub}^{s}=4$
are plotted in~\cref{fig:curve1a}. The approximated boundary is continuous and smooth
except in places with high curvature.~\Cref{fig:curve1b} shows that this overall
choice of parameters is sufficient to apply the constraint along the point cloud.

\begin{figure}[!htbp]
  \centering
  \begin{subfigure}[b]{0.48\textwidth}
    \centering
    \includegraphics[scale=1.0]{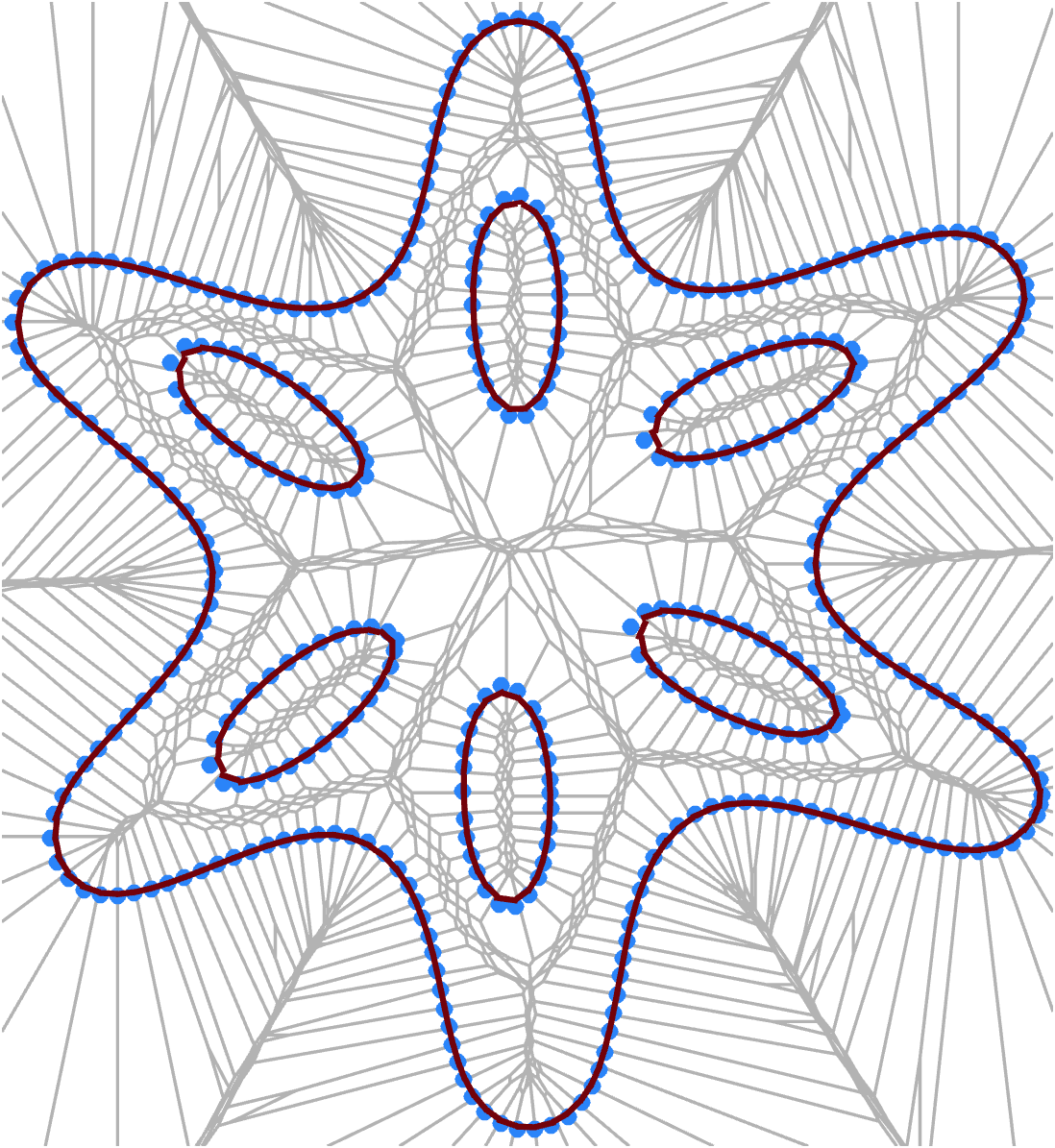}
    \caption{Sample points (blue dots), Voronoi diagram for k=4 (gray lines) and
      $\Gamma_{D,\text{sharp}}$ (red lines)}
    \label{fig:curve1a}
  \end{subfigure}
  \hfill
  \begin{subfigure}[b]{0.48\textwidth}
    \centering
    \includegraphics[scale=1.0]{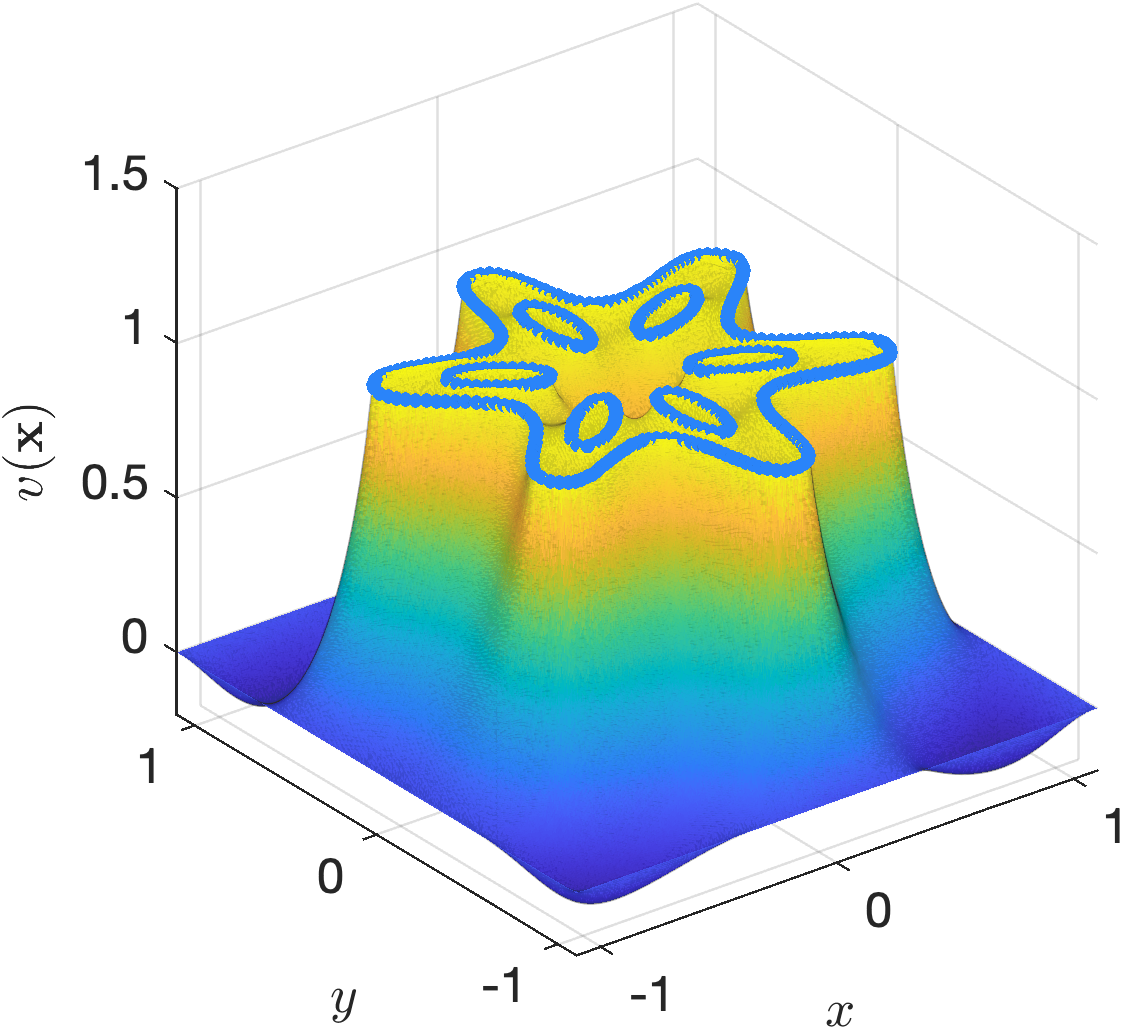}
    \caption{Solution to membrane structure problem}
    \label{fig:curve1b}
  \end{subfigure}
  \caption{\revb{Example \texttt{mc4.txt}}}
  \label{fig:curve1}
\end{figure}

The second example \texttt{oc16.txt} is shown in~\cref{fig:curve2a}. It consists of one 'snakelike' open curve with many non-smooth kinks sampled with slightly varying sampling density. Considering \cref{sec:robLim}
the characteristics of the approximated boundary on the kinks and open ends are of
interest here.

The approximated boundary for $\boldsymbol{r}=0.035$, $n_{query}^{s}=5$,
$l_{max}^{s}=8 \cdot 10^{-2}$ and $n_{sub}^{s}=4$ is shown in~\cref{fig:curve2b}.
Separating the two groups of Voronoi regions that likewise exist in this case is not
possible. We opt to choose $\boldsymbol{r}$ relatively tight to filter all non
wanted Voronoi regions in the vicinity of the medial axis. Subsequently we miss one
Voronoi region that should be included (green square in~\cref{fig:curve2b}). While
the open ends and kinks are handled satisfactory the approximated boundary shows
several small discontinuities. 

Clearly the sharp interface approach does not in general yield watertight, continuous
reconstructions as one might wish for in a e.g. in a computer graphics context. But
this was also not our original intent. Our goal was to apply Dirichlet boundary
conditions on domains defined by point clouds in the context of partial differential
equations. The implied smoothness of possible solutions to such equation systems in
combination with the weak nature of the penalty method allow for reduced requirements
for the approximated boundary. 

\begin{figure}[!htbp]
  \centering
  \begin{subfigure}[b]{0.48\textwidth}
    \centering
    \includegraphics[scale=1.0]{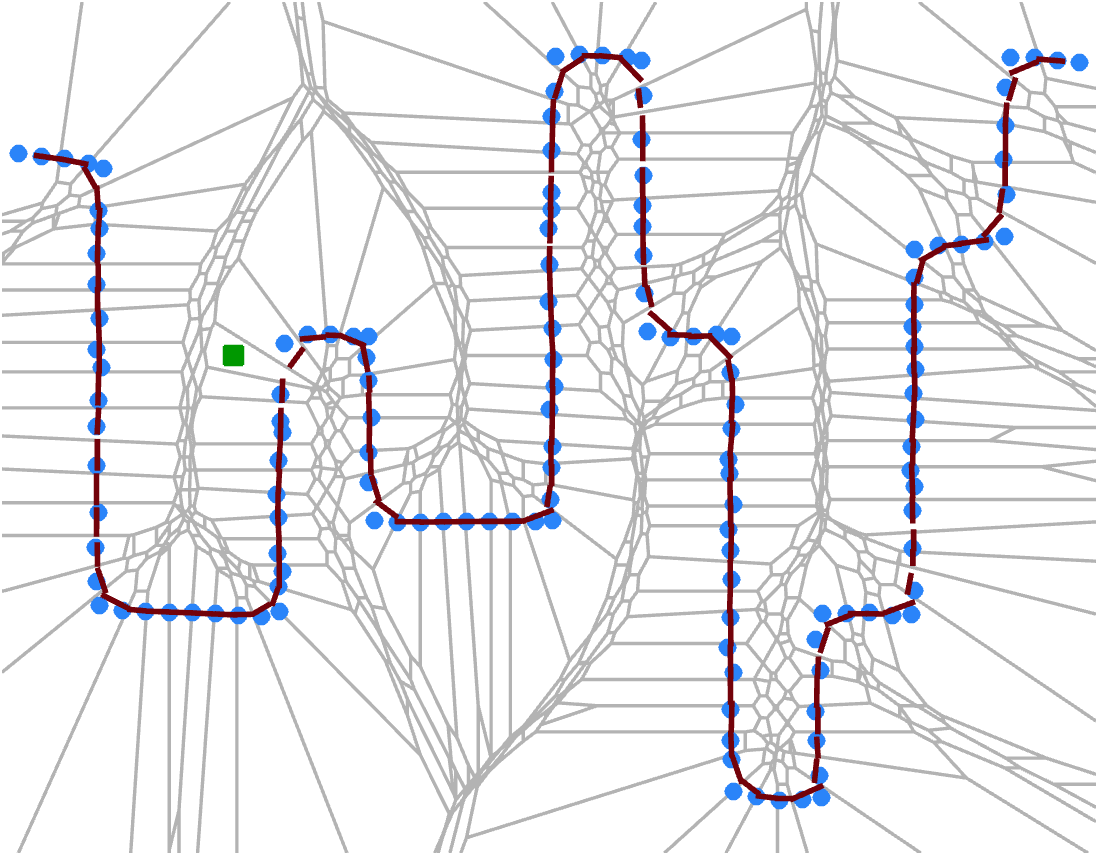}
    \caption{Sample points (blue dots), Voronoi diagram for k=4 (gray lines) and
      $\Gamma_{D,\text{sharp}}$ (red lines)}
    \label{fig:curve2a}
  \end{subfigure}
  \hfill
  \begin{subfigure}[b]{0.48\textwidth}
    \centering
    \includegraphics[scale=1.0]{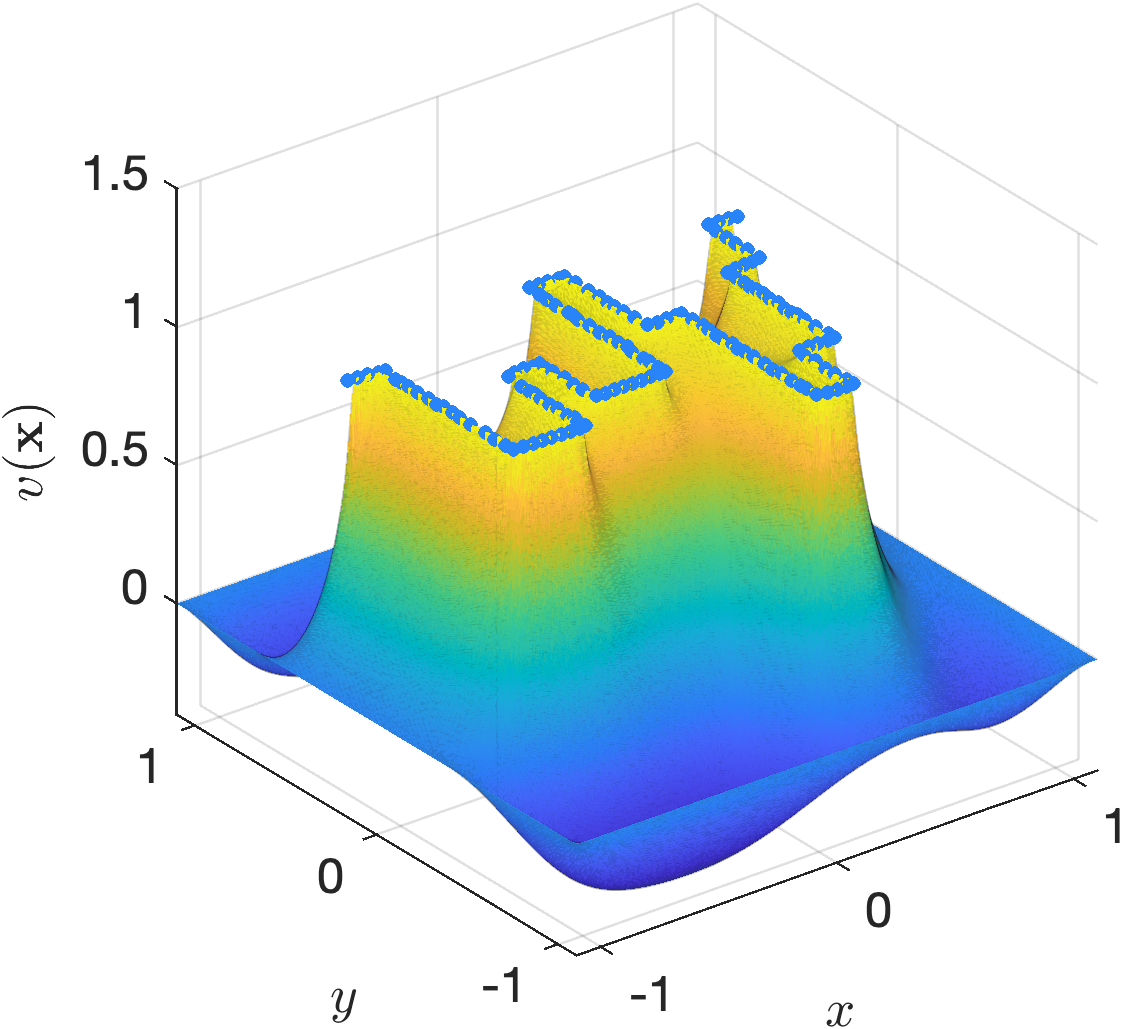}
    \caption{Solution to membrane structure problem}
    \label{fig:curve2b}
  \end{subfigure}
  \caption{\revb{Example \texttt{oc16.txt}}}
  \label{fig:curve2}
\end{figure}

\emph{Remark:} Remember that the Voronoi diagram is actually never
constructed. It is only approximated implicitly from $N^{k}(P, \boldsymbol{x})$ as
described by~\cref{fig:sharpInt1} and~\cref{fig:sharpInt2}. Nevertheless, we plot
its outline in gray because this gives a good impression on how the red sharp
boundary boundary is reconstructed.}


\section{Conclusions} \label{sec:conclusions}

Solving contour integrals on non-boundary-conforming meshes is essential to apply
boundary conditions in embedded domain methods such as the finite cell method. In
\cite{Kudela2020} an approach is presented where the domain is solely defined by a set of points. Therein, Neumann boundary conditions are applied by implicitly approximating non-zero Neumann boundaries using the diffuse interface approach. 

This contribution builds upon this idea but modifies it in such a way that the same local approximation can be
integrated without assuming a diffuse interface. In essence this can be achieved by
closely examining the implicit boundary approximation and interpreting it with
respect to order-$k$ Voronoi diagrams. Both approaches are compared in a numerical
example for the case of enforcing Dirichlet boundary conditions with the penalty
method. Therein it becomes apparent that the diffuse interface approach inherently suffers from
a trade-off between integration effort and accuracy which is completely eliminated in the newly presented sharp interface approach. \reva{Therefore, the sharp interface approach should be preferred wherever possible.}


\newpage
\section*{Acknowledgements} 
The last author gratefully acknowledges the financial support of the German Research Foundation (DFG) under
grant KO 4570/1-1.

\bibliographystyle{ieeetr}
\bibliography{zot_references,zot_references_stefan}

\begin{thebibliography}{10}

\bibitem{Mantyla1988}
M.~M{\"a}ntyl{\"a}, {\em An introduction to solid modeling}.
\newblock No.~13 in Principles of computer science series, {Rockville}:
  {Computer Science Press}, 1988.

\bibitem{Wassermann2019}
B.~Wassermann, S.~Kollmannsberger, S.~Yin, L.~Kudela, and E.~Rank,
  ``Integrating {{CAD}} and numerical analysis: `{{Dirty}} geometry' handling
  using the {{Finite Cell Method}},'' {\em Computer Methods in Applied
  Mechanics and Engineering}, vol.~351, pp.~808--835, July 2019.

\bibitem{Massarwi2016}
F.~Massarwi and G.~Elber, ``A {{B}}-spline based framework for volumetric
  object modeling,'' {\em Computer-Aided Design}, vol.~78, pp.~36--47, Sept.
  2016.

\bibitem{Zhang2016a}
Y.~J. Zhang, {\em Geometric {{Modeling}} and {{Mesh Generation From Scanned
  Images}}}.
\newblock Chapman \& {{Hall}}/{{CRC}} mathematical and computational imaging
  sciences, {Chapman \& Hall/CRC}, 2016.

\bibitem{deGeus2017}
T.~{de Geus}, J.~Vond{\v r}ejc, J.~Zeman, R.~Peerlings, and M.~Geers, ``Finite
  strain {{FFT}}-based non-linear solvers made simple,'' {\em Computer Methods
  in Applied Mechanics and Engineering}, vol.~318, pp.~412--430, May 2017.

\bibitem{Korshunova2020}
N.~Korshunova, J.~Jomo, G.~L{\'e}k{\'o}, D.~Reznik, P.~Bal{\'a}zs, and
  S.~Kollmannsberger, ``Image-based material characterization of complex
  microarchitectured additively manufactured structures,'' {\em Computers \&
  Mathematics with Applications}, vol.~80, pp.~2462--2480, Dec. 2020.

\bibitem{Korshunova2021a}
N.~Korshunova, G.~Alaimo, S.~Hosseini, M.~Carraturo, A.~Reali, J.~Niiranen,
  F.~Auricchio, E.~Rank, and S.~Kollmannsberger, ``Bending behavior of
  octet-truss lattice structures: {{Modelling}} options, numerical
  characterization and experimental validation,'' {\em Materials \& Design},
  vol.~205, p.~109693, July 2021.

\bibitem{Gravenkamp2020}
H.~Gravenkamp, A.~A. Saputra, and S.~Eisentr{\"a}ger, ``Three-dimensional
  image-based modeling by combining {{SBFEM}} and transfinite element shape
  functions,'' {\em Computational Mechanics}, vol.~66, pp.~911--930, Oct. 2020.

\bibitem{Cottrell2009}
J.~A. Cottrell, T.~J.~R. Hughes, and Y.~Bazilevs, {\em Isogeometric
  {{Analysis}}: {{Toward Integration}} of {{CAD}} and {{FEA}}}.
\newblock {John Wiley \& Sons}, Aug. 2009.

\bibitem{Kudela2020}
L.~Kudela, S.~Kollmannsberger, U.~Almac, and E.~Rank, ``Direct structural
  analysis of domains defined by point clouds,'' {\em Computer Methods in
  Applied Mechanics and Engineering}, vol.~358, p.~112581, Jan. 2020.

\bibitem{Saulev1963}
V.~Saul'ev, ``{Saul'ev, V.K.: On solution of some boundary value problems on
  high performance computers by fictitious domain method. Sibirian Math. J. 4,
  912\textendash 925 (1963)},'' {\em Sibirian Mathematical Journal}, vol.~4,
  pp.~912--925, 1963.

\bibitem{Burman2015b}
E.~Burman, S.~Claus, P.~Hansbo, M.~G. Larson, and A.~Massing, ``{{CutFEM}}:
  {{Discretizing}} geometry and partial differential equations,'' {\em
  International Journal for Numerical Methods in Engineering}, vol.~104,
  pp.~472--501, Nov. 2015.

\bibitem{Breitenberger2015}
M.~Breitenberger, A.~Apostolatos, B.~Philipp, R.~W{\"u}chner, and K.~U.
  Bletzinger, ``Analysis in computer aided design: {{Nonlinear}} isogeometric
  {{B}}-{{Rep}} analysis of shell structures,'' {\em Computer Methods in
  Applied Mechanics and Engineering}, vol.~284, pp.~401--457, Feb. 2015.

\bibitem{Kamensky2015}
D.~Kamensky, M.-C. Hsu, D.~Schillinger, J.~A. Evans, A.~Aggarwal, Y.~Bazilevs,
  M.~S. Sacks, and T.~J.~R. Hughes, ``An immersogeometric variational framework
  for fluid\textendash structure interaction: {{Application}} to bioprosthetic
  heart valves,'' {\em Computer Methods in Applied Mechanics and Engineering},
  vol.~284, pp.~1005--1053, Feb. 2015.

\bibitem{Badia2018}
S.~Badia, F.~Verdugo, and A.~F. Mart{\'i}n, ``The aggregated unfitted finite
  element method for elliptic problems,'' {\em Computer Methods in Applied
  Mechanics and Engineering}, vol.~336, pp.~533--553, 2018.

\bibitem{Duster2008}
A.~D{\"u}ster, J.~Parvizian, Z.~Yang, and E.~Rank, ``The finite cell method for
  three-dimensional problems of solid mechanics,'' {\em Computer Methods in
  Applied Mechanics and Engineering}, vol.~197, pp.~3768--3782, Aug. 2008.

\bibitem{Duster2017}
A.~D{\"u}ster, E.~Rank, and B.~A. Szab{\'o}, ``The p-version of the finite
  element method and finite cell methods,'' in {\em Encyclopedia of
  {{Computational}} mechanics} (E.~Stein, R.~Borst, and T.~J.~R. Hughes, eds.),
  vol.~2, pp.~1--35, {Chichester, West Sussex}: {John Wiley \& Sons}, 2017.

\bibitem{Longva2020}
A.~Longva, F.~L{\"o}schner, T.~Kugelstadt, J.~A. {Fern{\'a}ndez-Fern{\'a}ndez},
  and J.~Bender, ``Higher-order finite elements for embedded simulation,'' {\em
  ACM Transactions on Graphics}, vol.~39, pp.~181:1--181:14, Nov. 2020.

\bibitem{babusFiniteElement}
I.~Babuska, ``The {{Finite Element Method}} with {{Penalty}},''

\bibitem{Schillinger2016a}
D.~Schillinger, I.~Harari, M.-C. Hsu, D.~Kamensky, S.~K.~F. Stoter, Y.~Yu, and
  Y.~Zhao, ``The non-symmetric {{Nitsche}} method for the parameter-free
  imposition of weak boundary and coupling conditions in immersed finite
  elements,'' {\em Computer Methods in Applied Mechanics and Engineering},
  vol.~309, pp.~625--652, Sept. 2016.

\bibitem{Kollmannsberger2015}
S.~Kollmannsberger, A.~{\"O}zcan, J.~Baiges, M.~Ruess, E.~Rank, and A.~Reali,
  ``Parameter-free, weak imposition of {{Dirichlet}} boundary conditions and
  coupling of trimmed and non-conforming patches,'' {\em International Journal
  for Numerical Methods in Engineering}, vol.~101, pp.~670--699, Mar. 2015.

\bibitem{amentNewVoronoibased}
N.~Amenta, M.~Bern, and M.~Kamvysselis, ``A new {{Voronoi}}-based surface
  reconstruction algorithm,'' in {\em Proceedings of the 25th Annual Conference
  on {{Computer}} Graphics and Interactive Techniques - {{SIGGRAPH}} '98},
  ({Not Known}), pp.~415--421, {ACM Press}, 1998.

\bibitem{liSOLVINGPDES}
X.~Li, J.~Lowengrub, A.~R. Tz, and A.~Voigt, ``{{SOLVING PDES IN COMPLEX
  GEOMETRIES}}: A {{DIFFUSE DOMAIN APPROACH}},'' p.~27.

\bibitem{ratzPDESurfaces}
A.~R{\"a}tz and A.~Voigt, ``{{PDE}}'s on surfaces---a diffuse interface
  approach,'' {\em Communications in Mathematical Sciences}, vol.~4, no.~3,
  pp.~575--590, 2006.

\bibitem{leeRegularizedDirac}
H.~G. Lee and J.~Kim, ``Regularized {{Dirac}} delta functions for phase field
  models,'' {\em International Journal for Numerical Methods in Engineering},
  vol.~91, pp.~269--288, July 2012.

\bibitem{aurenVoronoiDiagrams-2}
F.~Aurenhammer, R.~Klein, and D.-T. Lee, {\em Voronoi Diagrams and {{Delaunay}}
  Triangulations}.
\newblock 2013.

\bibitem{schmiOrderkVoronoi}
D.~Schmitt and J.-C. Spehner, ``Order-k {{Voronoi Diagrams}}, k-{{Sections}},
  and k-{{Sets}},'' in {\em Discrete and {{Computational Geometry}}} (G.~Goos,
  J.~Hartmanis, J.~{van Leeuwen}, J.~Akiyama, M.~Kano, and M.~Urabe, eds.),
  vol.~1763, pp.~290--304, {Berlin, Heidelberg}: {Springer Berlin Heidelberg},
  2000.

\bibitem{deyVoronoibasedFeature}
T.~K. Dey and L.~Wang, ``Voronoi-based feature curves extraction for sampled
  singular surfaces,'' {\em Computers \& Graphics}, vol.~37, pp.~659--668, Oct.
  2013.

\bibitem{edelsShapeSet}
H.~Edelsbrunner, D.~Kirkpatrick, and R.~Seidel, ``On the shape of a set of
  points in the plane,'' {\em IEEE Transactions on Information Theory},
  vol.~29, pp.~551--559, July 1983.

\bibitem{allieVoronoibasedVariational}
P.~Alliez, D.~{Cohen-Steiner}, Y.~Tong, and M.~Desbrun, ``Voronoi-based
  {{Variational Reconstruction}} of {{Unoriented Point Sets}},'' p.~10.

\bibitem{Ohrhallinger2021}
S.~Ohrhallinger, J.~Peethambaran, A.~D. Parakkat, T.~K. Dey, and
  R.~Muthuganapathy, ``{{2D Points Curve Reconstruction Survey}} and
  {{Benchmark}},'' {\em Computer Graphics Forum}, vol.~40, pp.~611--632, May
  2021.

\bibitem{Zander2014}
N.~Zander, T.~Bog, M.~Elhaddad, R.~Espinoza, H.~Hu, A.~Joly, C.~Wu, P.~Zerbe,
  A.~D{\"u}ster, S.~Kollmannsberger, J.~Parvizian, M.~Ruess, D.~Schillinger,
  and E.~Rank, ``{{FCMLab}}: {{A}} finite cell research toolbox for
  {{MATLAB}},'' {\em Advances in Engineering Software}, vol.~74, pp.~49--63,
  Aug. 2014.

\bibitem{zandeFiniteCell}
N.~Zander, S.~Kollmannsberger, M.~Ruess, Z.~Yosibash, and E.~Rank, ``The
  {{Finite Cell Method}} for linear thermoelasticity,'' {\em Computers \&
  Mathematics with Applications}, vol.~64, pp.~3527--3541, Dec. 2012.

\bibitem{nguyeDiffuseNitsche}
L.~H. Nguyen, S.~K. Stoter, M.~Ruess, M.~A. Sanchez~Uribe, and D.~Schillinger,
  ``The diffuse {{Nitsche}} method: Dirichlet constraints on phase-field
  boundaries,'' {\em International Journal for Numerical Methods in
  Engineering}, vol.~113, pp.~601--633, Jan. 2018.

\bibitem{Zander2015}
N.~Zander, T.~Bog, S.~Kollmannsberger, D.~Schillinger, and E.~Rank,
  ``Multi-level hp-adaptivity: high-order mesh adaptivity without the
  difficulties of constraining hanging nodes,'' {\em Computational Mechanics},
  vol.~55, pp.~499--517, Feb. 2015.

\bibitem{DAngella2018}
D.~D'Angella, S.~Kollmannsberger, E.~Rank, and A.~Reali, ``Multi-level
  {{B\'ezier}} extraction for hierarchical local refinement of {{Isogeometric
  Analysis}},'' {\em Computer Methods in Applied Mechanics and Engineering},
  vol.~328, pp.~147--174, Jan. 2018.

\bibitem{Kopp2021a}
P.~Kopp, E.~Rank, V.~M. Calo, and S.~Kollmannsberger, ``Efficient multi-level
  hp-finite elements in arbitrary dimensions,'' {\em arXiv:2106.08214 [cs,
  math]}, June 2021.

\bibitem{burgeAnalysisDiffuse}
M.~Burger, O.~L. Elvetun, and M.~Schlottbom, ``Analysis of the {{Diffuse Domain
  Method}} for {{Second Order Elliptic Boundary Value Problems}},'' {\em
  Foundations of Computational Mathematics}, vol.~17, pp.~627--674, June 2017.

\end{thebibliography}

\end{document}